\documentclass[12pt,a4paper]{report}
\input{diagrams}
\diagramstyle[height=1.5em,scriptlabels,PS=dvips]

\usepackage{amsmath}
\usepackage{amssymb}
\usepackage{theorem,french}
\usepackage{xy}

\xyoption{all} 
\CompileMatrices 

\newarrow{To}----{->}
\newarrow{Into}C---{->}
\newarrow{Suto}----{->>}
\newarrow{impli}===={=>}

\font\tencyr=wncyr10
\font\sevencyr=wncyr7
\font\fivecyr=wncyr5

\newfam\cyrfam
\textfont\cyrfam\tencyr
\scriptfont\cyrfam\sevencyr
\scriptscriptfont\cyrfam\fivecyr

\def\gp{\mathfrak{p}}
\def\gt{\mathfrak{t}}

\def\gd{\mathfrak{d}}
\def\gr{\mathfrak{r}}
\def\ga{\mathfrak{a}}

\def\N{{\bf N}}
\def\M{\text{\bf M}}
\def\Z{{\bf Z}}
\def\H{{\bf H}}
\def\Q{{\bf Q}}
\def\P{\Bbb P}
\def\A{\Bbb A}

\def\F{\text{\bf F}}

\def\C{\mbox{\bf C}}

\def\Gl{{\rm Gl}}

\newcommand{\eps}{\varepsilon}

\newcommand{\Ker}{\hbox{{\rm ker}}}

\newcommand{\Ens}{\hbox{{\rm Ens}}}

\newcommand{\tors}{\hbox{{\rm tors}}}

\newcommand{\Id}{\text{{\rm Id}}}

\newcommand{\Spec}{\hbox{{\rm Spec}}}
\newcommand{\eme}{\text{\rm i\`eme}}
\newcommand{\res}{\hbox{{\rm res}}}

\newcommand{\PGL}{\hbox{{\rm PGL}}}

\newcommand{\Div}{\hbox{{\rm Div}}}
\newcommand{\Hom}{\hbox{{\rm Hom}}}
\newcommand{\Fix}{\hbox{{\rm Fix}}}
\newcommand{\Sch}{\hbox{{\rm Sch}}}
\newcommand{\Pic}{\hbox{{\rm Pic}}}
\newcommand{\Supp}{\hbox{{\rm Supp}}}
\newcommand{\glo}{\text{\rm gl}}

\newcommand{\loc}{\text{\rm loc}}

\newcommand{\Krull}{\text{\rm Krull}}

\newcommand{\Aut}{\hbox{{\rm Aut}}}

\newcommand{\im}{\mathopen{\mathchoice
             {\hbox{{\rm Im}}}
             {\hbox{{\rm Im}}}
             {\hbox{\small {\rm Im}}}
             {\hbox{\tiny {\rm Im}}}
             }}

\def\theenumi{(\roman{enumi})}

\def\theenumii{(\alph{enumii})}
\def\p@enumii{\theenumi}

\def\theenumiii{\Alph{enumiii}}
\def\p@enumiii{\theenumi\theenumii}

\def\p@enumiv{\p@enumiii\theenumiii.}

\newif\ifnormalesBeweisEnde

  {\vskip 0.3ex plus 0.5ex minus 0ex \pagebreak[1]
   \global\normalesBeweisEndetrue
   \trivlist
   \item[\hskip\labelsep {\textsc{Proof} \rm of #1}:]}%
  {\ifnormalesBeweisEnde \EndOfBeweis \fi
   \endtrivlist
   \vskip 1ex plus 1ex minus 0ex \pagebreak[2]}
\newenvironment{pr}
  {\vskip 0.3ex plus 0.5ex minus 0ex \pagebreak[1]
   \global\normalesBeweisEndetrue
   \trivlist
   \item[\hskip\labelsep \textsc{Preuve} :]}%
  {\ifnormalesBeweisEnde \EndOfBeweis \fi
   \endtrivlist
   \vskip 1ex plus 1ex minus 0ex \pagebreak[2]}
  {\ifnormalesBeweisEnde \EndOfBeweis \fi
   \endtrivlist
   \vskip 1ex plus 1ex minus 0ex \pagebreak[2]}
\def\EndOfBeweis{\hskip .5em \vrule width 1.0ex height 1.0ex depth 0.3ex}

\advance\topmargin by-1cm 
\advance\textheight by 0.25in 
\advance\oddsidemargin by-1.2cm 
\advance\evensidemargin by-1.2cm 
\advance\textwidth by 2.4cm 
\theoremstyle{plain} 
\newtheorem{thma}{Theorem}[subsection]
{\theorembodyfont{\rmfamily}
\newtheorem{dea}[thma]{Definition}

\newtheorem{ex}[thma]{Exemple}

\newtheorem{rem}[thma]{Remarque}}
\newtheorem{prop}[thma]{Proposition}

\newtheorem{lem}[thma]{Lemme}
\newtheorem{thm}[thma]{Th\'eor\`eme}

\newtheorem{cor}[thma]{Corollaire}

\newtheorem{thmfs}{Theorem}[section]
{\theorembodyfont{\rmfamily}
\newtheorem{thmf}[thmfs]{Th\'eor\`eme}
\newtheorem{lemf}[thmfs]{Lemme}
\newtheorem{propf}[thmfs]{Proposition}

}

{\obeylines
\gdef\iffin{\parskip0pt\parindent0pt 
            \vskip1cm
            \footnotesize
            courriel : jose.bertin@ujf-grenoble.fr, ariane.mezard@ujf-grenoble.fr
            \bigskip
            Universit\'e de Grenoble I
            {\bf Institut Fourier}
            UMR 5582 CNRS-UJF
            UFR de Math\'ematiques
            B.P. 74
            38402 Saint-Martin d'H\`eres Cedex (France)
            }}
\begin{document}
/usr/local/tex/lib/texmf/tex/latex/french/inputs/french.sty

\Large
$$\mbox{D\'eformations formelles des rev\^etements sauvagement}$$
$$\mbox{ramifi\'es de courbes alg\'ebriques}$$
\large
$$\mbox{Jos\'e Bertin et Ariane M\'ezard}$$\\
\noindent
\normalsize
{\def\thefootnote{\relax}\footnote{{\bf Mots cl\'es } : anneau de d\'eformation versel,
  obstruction, rev\^etement de courbes alg\'ebriques,
    ramification sauvage, \'equations d'Artin-Schreier, espaces de
    modules.}\addtocounter{footnote}{-1}}
{\bf Abstract}: In this paper we study formal moduli for wildly
ramified Galois covering. We prove a local-global principle.
We then focus on the infinitesimal deformations of the $\Z/p\Z$-covers.
We explicitly compute a deformation of an automorphism of order $p$
which implies a universal obstruction for $p>2$. 
By deforming Artin-Schreier equations we obtain a lower bound on
the dimension of the local versal deformation ring.
At last, by
comparing the global versal deformation ring to the complete local ring
in a point of a moduli space, we determine the dimensions of the
global and local versal deformation ring.
\section{Introduction}
Soit $k$ un corps alg\'ebriquement clos de caract\'eristique $p>0$ et
soit $W(k)$ l'anneau des vecteurs de Witt de $k$. Soit
$\pi:C\rightarrow\Sigma$ un rev\^etement d\'efini sur $k$ entre deux courbes
alg\'ebriques compl\`etes et lisses. Nous cherchons \`a relever ce
rev\^etement \`a un anneau de valuation discr\`ete complet $R$
dominant $W(k)$. Des motivations pour aborder un tel probl\`emes
peuvent \^etre trouv\'ees dans l'expos\'e de Oort (\cite{Oo}) et dans
les travaux r\'ecents de Green et Matignon (\cite{GrMa1}, \cite{GrMa2}).\\
Nous nous limitons au cas o\`u le rev\^etement est galoisien de groupe
$G$. Le rel\`evement de $\pi$ \`a $R$ est \'equivalent au
rel\`evement \`a $R$ de la courbe $C$ et de l'action du groupe $G$. Lorsque
les groupes d'inertie sont d'ordre premier \`a $p$, autrement dit si
la ramification est mod\'er\'ee, nous savons que le probl\`eme de
rel\`evement a une solution (\cite{Gro3}, \cite{Oo}). En revanche
s'il y a de la ramification sauvage, ce probl\`eme peut conduire \`a
des obstructions.\\
Ces obstructions sont de nature locales (voir notamment \cite{Be}, \cite{GrMa1}). Par
exemple elles peuvent porter sur le caract\`ere d'Artin attach\'e \`a
un point $x$ de $C$ dont le groupe d'inertie $G_x$ a un ordre divisible par
$p$. En g\'en\'eral, elles s'expriment comme obstructions \`a la
d\'eformation d'un sous-groupe fini du groupe des automorphismes de
$k[[T]]$. Sekiguchi, Oort et Suwa
(\cite{SeOoSu}) ont \'etudi\'e le probl\`eme du rel\`evement de $\pi$ en travaillant sur les
d\'eformations globales de $C$ et $G$. Nous adoptons un autre point de
vue : nous introduisons les
d\'eformations locales pour lesquelles la courbe $C$ est remplac\'ee par un voisinage
formel d'un point fixe et $G$ par le groupe d'inertie $G_x$ (d'ordre
divisible par $p$) correspondant.\\ 
Notre premier r\'esultat (Th\'eor\`eme \ref{thmloglo}) est l'\'enonc\'e d'un
principe local-global qui relie les d\'eformations globales des
d\'eformations locales. Un r\'esultat de cette nature avait d\'eja
\'et\'e utilis\'e par Green et Matignon (\cite{GrMa1}).
D'apr\`es la th\'eorie classique des d\'eformations, les
d\'eformations locales (ou globales) sont d\'ecrites par un anneau
local complet. Notons $R_{G_x}$ l'anneau associ\'e aux d\'eformations
locales. Gr\^ace au principe local-global, nous ramenons l'\'etude du
rel\`evement de $\pi$ \`a l'\'etude de l'anneau $R_{G_x}$. Nous
pouvons rapprocher le probl\`eme de la description de $R_{G_x}$ \`a
celui introduit par Mazur portant sur les anneaux de d\'eformations universels
de repr\'esentations galoisiennes (\cite{Ma},\cite{Me}). Cette
description est en g\'en\'eral compliqu\'ee.
Nous identifions l'espace
tangent de $R_{G_x}$ au groupe $H^1(G_x,\Theta)$, avec $\Theta$ le $G_x$-module des champs de vecteurs formels au
point $x$. Puis nous montrons que les obstructions appartiennent \`a $H^2(G_x,\Theta)$. Lorsque $G_x$ est un $p$-groupe cyclique, nous
calculons les dimensions de ces groupes de cohomologies.\\
Notre second r\'esultat porte sur la description de $R_{G_x}$ pour $G_x$
cyclique d'ordre $p$. Le rel\`evement en caract\'eristique
z\'ero est alors possible (\cite{GrMa1},\cite{SeOoSu}). Notons $m$ le conducteur de Hasse de $G_x$. En mettant en \'evidence une obstruction cohomologique non triviale
universelle, nous prouvons que si $p>2$ et si $(m,p)\not=(1,3)$, alors $R_{G_x}$ est toujours singulier. 
Par un argument de d\'eformations des \'equations
d'Artin-Schreier, nous minorons la dimension de Krull
de $R_{G_x}$. Le principe local-global et un r\'esultat de
Harbater (\cite{Har}) permettent alors de prouver que cette borne
inf\'erieure est effectivement la dimension de Krull. Si le conducteur
est $m=1$, l'anneau $R_{G_x}$ est d\'ecrit explicitement \`a l'aide d'une
unique \'equation que nous identifions \`a un polyn\^ome de Tchebychev.
Enfin pour les petites valeurs de $m$ ($m<p-1$), nous observons que
$R_{G_x}$ est d'intersection compl\`ete.
\section{D\'eformations locales et globales}
Soit $\Lambda$ un anneau de valuation discr\`ete complet de
caract\'eristique 0 de corps r\'esiduel $k$ ; nous pouvons 
choisir $\Lambda=W(k)$.
Consid\'erons la cat\'egorie  $\widehat{\cal C}$ dont les objets sont les
$\Lambda$-alg\`ebres locales $R$ noetheriennes compl\`etes, avec $k\cong R/{\cal M}_R$ pour ${\cal M}_R$ id\'eal
maximal de $R$ ; les fl\`eches
sont les $\Lambda$-morphismes d'anneaux locaux.
La cat\'egorie ${\cal C}$ est la sous-cat\'egorie pleine de $\widehat{\cal C}$
form\'ee des objets  qui sont de longueur finie comme 
$\Lambda$-modules. En particulier notons $k[\eps]$, o\`u $\eps^2=0$, l'anneau de ${\cal C}$ des nombres duaux sur $k$.
Une surjection $A'\rightarrow A$ dans ${\cal C}$ est
une {\it petite extension} si son
noyau $\ga$ est principal et v\'erifie 
$\ga{\cal M}_{A'}=0$.\\

\noindent 
Soit $C\rightarrow C/G$ un {\it rev\^etement sauvagement
  ramifi\'e},\index{rev\^etement sauvagement ramifi\'e} 
que nous identifierons au couple $(C,G)$, compos\'e
d'une $k$-courbe alg\'ebrique projective lisse $C$ et d'un sous-groupe
fini $G$ de Aut$_k(C)$ d'ordre divisible par $p$.
Un point $x\in C$ est un {\it point de ramification sauvage}\index{point de ramification sauvage} si le
stabilisateur (le sous-groupe d'inertie) $G_x$ au point $x$ a un
ordre divisible par $p$. Ces points vont jouer le r\^ole de points
singuliers et vont contribuer aux d\'eformations du rev\^etement. \\

Rappelons bri\`evement la terminologie usuelle (\cite{DeMu},\cite{La},\cite{Sc}). 
Une {\it d\'eformation  de} $(C,G)$ \`a $A$ objet de ${\cal C}$ est
une classe d'isomorphismes $G$-\'equivariants induisant l'identit\'e sur
$C$ de couples $(X,G)$ d\'efinis par\\
-un {\it rel\`evement de} $C$ \`a $A$, i.e. une courbe propre et lisse
$\gp : X\rightarrow \mbox{Spec}A$ avec un isomorphisme $\gp^{-1}(0)\cong
C$ (c'est-\`a-dire la fibre de $X$ au-dessus du point ferm\'e de
Spec$A$ est $C$).\\
-un rel\`evement de $G$ \`a $X$ en un sous-groupe de Aut$_A(X)$, tel
que l'isomorphisme $C\cong \gp^{-1}(0)=X\otimes_A k$ soit 
$G$-\'equivariant.\\
Nous d\'efinissons ainsi un foncteur covariant
$$D_{\glo}:{\cal C}\rightarrow
\Ens,\;\;\;A\mapsto\{\mbox{d\'eformations de } (C,G) \mbox{ \`a } A\}$$ 
 Ce foncteur est dit {\it foncteur des d\'eformations
  de} $(C,G)$.\index{foncteur de d\'eformations} L'indice gl de $D_{\glo}$ indique que nous d\'eformons
globalement la courbe.\\
Soit deux foncteurs covariants $D,D':{\cal C}\rightarrow \Ens$.
Un morphisme de foncteurs $\phi:D\rightarrow D'$ est {\it lisse} si
pour toute petite extension de ${\cal C}$ $A\rightarrow A'$, $\phi$ induit un morphisme
surjectif
$$D'(A')\twoheadrightarrow D(A')\times_{D(A)}D'(A)$$
 Le foncteur covariant $D$ 
{\it admet une d\'eformation
 verselle}\index{d\'eformation ! verselle} si : (\cite{Sc} \S 2)
Il existe $R$ objet de $\widehat{\cal C}$ tel que
$D(k[\eps])\cong {\cal M}_R/{\cal M}_R^2$,  et un morphisme lisse de foncteurs
$\xi:D\rightarrow \Hom_{\Lambda}(R,\cdot)$.
Le couple $(R,\xi)$ est dit {\it
  d\'eformation
verselle} et l'anneau $R$ est dit {\it anneau de d\'eformations
versel}.\index{anneau de d\'eformations ! versel}
Si $\xi$ est un isomorphisme de foncteur, alors 
$(R,\xi)$ est dit
{\it d\'eformation universelle},\index{d\'eformation ! universelle} et
l'anneau $R$ est dit {\it anneau de
  d\'eformations universel}\index{anneau de d\'eformations ! universel}.\\
Soit ${\cal T}_{C}$ le faisceau tangent \`a la courbe, dual du faisceau des diff\'erentielles $\Omega^1_{C/k}$.
  Comme ${\cal T}_C$ est un $G$-faisceau, les groupes de cohomologie
  $H^i(G,{\cal T}_C)$ sont naturellement des $G$-modules.
\begin{thmf}
\label{threpre}
Le foncteur des d\'eformations $D_{\glo}$ de $(C,G)$ admet une
d\'eformation verselle, et si $H^0(C,{\cal T}_{C})^G=0$, cette d\'eformation
est universelle. 
\end{thmf}
\begin{pr}
Il suffit de reprendre la d\'emonstration de Schlessinger (\cite{Sc} \S 3.7)
qui \'etablit l'existence d'une d\'eformation verselle du foncteur des
d\'eformations de $C$, en constatant que les constructions faites sont
compatibles avec l'action de $G$.
\end{pr}
L'hypoth\`ese du th\'eor\`eme \ref{threpre} est v\'erifi\'ee par
exemple si la courbe $C$ est de genre $g\geq2$, car alors $H^0(C,{\cal T}_C)=0$
et dans ce cas $D_{\glo}$ admet une d\'eformation universelle. Par le th\'eor\`eme d'alg\'ebrisation de
Grothendieck (FGA, \cite{Gr2}) le sch\'ema formel qui pro-repr\'esente
$D_{\glo}$ est alg\'ebrisable, donc est le compl\'et\'e formel d'une courbe
propre et lisse sur $\Spec R$ uniquement d\'efinie. De plus l'action
de $G$ provient d'une action sur cette courbe.
Nous allons \`a pr\'esent d\'efinir le foncteur des d\'eformations infinit\'esimales.
Soit $y\in C/G$ image d'un
point de ramification sauvage $x\in C$. Notons $G_x$ le stabilisateur
de $x$, $\widehat{\cal O}_{C,x}\cong k[[T]]$ l'anneau local
compl\'et\'e de la courbe
$C$ au point $x$ et
$\widehat{\cal T}_{C,x}$ la fibre compl\'et\'ee du faisceau tangent au point $x$. Il en d\'ecoule une ``repr\'esentation'' de $G_x$ dans $k[[T]]$, c'est-\`a-dire un morphisme injectif $G_x\hookrightarrow\mbox{Aut}k[[T]]$. Un repr\'esentant
$(X,G)$ d'une d\'eformation de $(C,G)$ \`a $A$ d\'efinit de m\^eme une
repr\'esentation $G_x\hookrightarrow \Aut A[[T]]$ qui par la
r\'eduction par le
morphisme canonique $A\rightarrow k\cong A/{\cal M}_A$ redonne la
repr\'esentation initiale de $G_x$ sur $k[[T]]$.\\
Ceci sugg\`ere de formuler le probl\`eme suivant :
fixons un sous-groupe fini de $\Aut k[[T]]$ not\'e $G$ en absence
d'ambigu\"\i t\'e avec le sous-groupe des automorphismes de la courbe $C$. 
Nous obtenons une repr\'esentation 
$\bar{\rho}:G\hookrightarrow\Aut k[[T]]$.
Pour $A$ objet de ${\cal C}$, d\'efinissons
le {\it foncteur des d\'eformations
  infinit\'esimales} $D_G$ de
$\bar{\rho}:G\hookrightarrow\mbox{Aut}k[[T]]$ :
$$D_G:{\cal C}\rightarrow \Ens\;\;\; A\mapsto\left\{\begin{array}{ll}
\mbox{rel\`evement } G\rightarrow \Aut A[[T]] \mbox{ modulo}\cr
\mbox{la conjugaison par un \'el\'ement}\cr
\mbox{de} \ker(\Aut A[[T]]\rightarrow \Aut k[[T]])\end{array}\right\}$$
Les propri\'et\'es n\'ecessaires \`a
l'application du crit\`ere de Schlessinger se v\'erifient
ais\'ement.
\begin{thmf}
Le foncteur $D_G$ des rel\`evements infinit\'esimaux de
$\bar{\rho}: G\hookrightarrow\Aut k[[T]]$ admet une d\'eformation
verselle.
\end{thmf}
Contrairement au cas global, ce rel\`evement versel n'est en g\'en\'eral pas universel.
\section{Principe local-global}
Dans ce paragraphe nous allons relier les d\'eformations globales du
rev\^etement induit par $(C,G)$ aux d\'eformations locales associ\'ees
aux points de branchement sauvage. Le r\'esultat obtenu peut \^etre
rapproch\'e du principe local-global
utilis\'e par Green et
Matignon \cite{GrMa1} dans le cadre de la g\'eom\'etrie rigide.
\subsection{Rappel sur la cohomologie \'equivariante}
Dans ce paragraphe, nous rappellons bri\`evement les notions de
cohomologie \'equivariante qui sont n\'ecessaires pour d\'ecrire en
termes cohomologiques le morphisme local-global ; pour plus de
d\'etails voir, par exemple \cite{Gr1}.\\

\noindent
Fixons $X$ un sch\'ema de type fini sur $k$ et
supposons que le groupe (fini) $G$ agit sur $X$ (les \'enonc\'es de
cette partie seront appliqu\'es \`a la courbe $C$).
Soit ${\cal O}_X$ le faisceau des fonctions sur $X$.
Supposons que tout point de $X$ admet un voisinage affine $G$-stable
(ce qui est en particulier v\'erifi\'e  si $X$ est projectif).
Alors le sch\'ema quotient $\Sigma=X/G$ est d\'efini ; d'apr\`es [KaMa] ce
passage au quotient est compatible aux changements de base.\\
Un $(G,{\cal O}_X)$-{\it module} est un ${\cal O}_X$-module
quasi-coh\'erent sur lequel $G$ agit. Soit $F$ un $(G,{\cal
  O}_X)$-module, nous avons une action naturelle de $G$ sur $H^q(X,F)$.
Si $\pi:X\rightarrow \Sigma=X/G$, nous pouvons consid\'erer le faisceau 
not\'e $\pi_*^G(F)$ sur $\Sigma$ d\'efini par :
si $V$ est un ouvert de $\Sigma$,
$$V\mapsto \Gamma(V,\pi_*(F))^G=\Gamma(\pi^{-1}(V),F)^G$$
Nous d\'efinissons deux foncteurs covariants et exacts \`a gauche sur
la cat\'egorie des $(G,{\cal O}_X)$-modules :
$$\pi_*^G \mbox{ et } \Gamma^G(X,\bullet)$$
avec $\Gamma^G(X,F)=\Gamma(X,F)^G $.
Les d\'eriv\'ees de ces foncteurs sont respectivement  un faisceau de
modules sur $X/G$ not\'e ${\cal X}^q(G,F)=R^q\pi_*^G(X,F)$ et un
groupe $H^q(G,F)=R^q\Gamma^G(X,F)$. Le $k$-espace vectoriel $H^q(G,F)$
est le  {\it groupe de cohomologie \'equivariante}\index{groupe de cohomologie \'equivariante} de $G$ \`a
coefficients dans le $G$-faisceau  $F$.
\begin{rem}{\rm
Vu que $\Gamma(\Sigma,\pi_*^G(F))=\Gamma(X,F)^G$, nous avons 
$\Gamma_X^G=\Gamma_{\Sigma}\circ\pi_*^G$ ; d'o\`u les deux suites 
spectrales convergentes d'aboutissement $H^{\bullet}(G,F) $ :
$$'E_2^{p,q}=H^p(\Sigma,R^q\pi_*^G(F))\Rightarrow H^{p+q}(G,F)$$
$$''E_2^{p,q}=H^p(G,H^q(X,F))\Rightarrow H^{p+q}(G,F)$$}
\end{rem}
Rappelons comment calculer la cohomologie $G$-\'equivariante par un proc\'ed\'e
\`a la \v{C}ech (voir \cite{Gr1}).
Soit ${\cal U}=\{U_i\}$ un recouvrement fini de $X$ par des ouverts affines $G$-stables.
Notons ${\cal C}^q({\cal U},F)$ le $q^{\eme}$-faisceau des germes de 
cocha\^\i nes de \v{C}ech ; par d\'efinition pour $U$ ouvert de X
$$\Gamma(U,{\cal C}^q({\cal
 U},F))=\prod_{\begin{array}{cc}\scriptstyle{(i_0,\cdots,i_q)}\\ 
\scriptstyle{ U_{i_0}\cap\cdots\cap U_{i_q}\not=\emptyset}\end{array}}\Gamma(U\cap U_{i_0}\cap\cdots\cap U_{i_q},F)$$
Rappellons que le complexe de
\v{C}ech ${\cal C}^q\stackrel \delta\rightarrow{\cal C}^{q+1}$ est une
r\'esolution de $F$,  o\`u la
diff\'erentielle de \v{C}ech est d\'efinie par :
$$\forall\omega\in{\cal C}^q,\;\;\;\;\;\;\delta(\omega)_{i_0\cdots
  i_{q+1}}=\sum_{\lambda=0}^{q+1}(-1)^{\lambda}\res_{U\cap U_{i_0}\cap
  \cdots \cap U_{i_{q+1}}}(\omega_{i_0\cdots \hat{i}_{\lambda}\cdots
  i_{q+1}})$$
Soit $\Gamma(X,{\cal C}^q)=C^q({\cal U},F)$ ;
c'est un $G$-module car ${\cal C}^q$ est un $G$-faisceau, et nous pouvons
  former la r\'esolution standard de ce $G$-module
$$C^{p,q}=C^p(G,C^q({\cal U},F))$$
La diff\'erentielle $d:C^{p,q}\rightarrow C^{p+1,q}$ de la cohomologie
des groupes et la diff\'erentielle $\delta:C^{p,q}\rightarrow
C^{p,q+1}$ font du groupe bigradu\'e $\{C^{p,q}\}$ un complexe double.
Notons 
$$\gd=d+\delta$$
la diff\'erentielle totale. Donc pour $\omega\in C^{p,q}$, nous avons
$$\gd(\omega)=d(\omega)+(-1)^p\delta(\omega)$$
Notons $(C^{\bullet},\gd)$ le complexe simple associ\'e au complexe double.
\begin{thm}
\label{thmcomp}
Le complexe $(C^{\bullet},D)$ permet de calculer la cohomologie
\'equivariante, et
$$H^q(C^{\bullet})\cong H^q(G,F)$$
\end{thm}
(Voir
\cite{Go},\cite{Gr1}). Gr\^ace au
th\'eor\`eme \ref{thmcomp}, nous pouvons d\'ecrire
les groupes de cohomologie $H^1(G,{\cal T}_C)$ et $H^2(G,{\cal T}_C)$.
Soit $\omega=\{\{\zeta_i^{\sigma}\},\{\delta_{ij}\}\}\in C^{1,0}\oplus
C^{0,1}$.
Nous avons  $\gd\omega=0$ si et seulement si 
$-\delta\{\zeta_i^{\sigma}\}+d\{\delta_{ij}\}=0$
et $\{\zeta_i^{\sigma}\}$ et
$\{\delta_{ij}\}$ sont des 1-cocycles :
$$\zeta_i^{\sigma\tau}=\zeta_i^{\sigma}+\sigma(\zeta_i^{\tau}),\;
\delta_{ik}=\delta_{ij}+\delta_{jk} \mbox{ sur }U_i\cap U_j\cap U_k$$
Donc
$\zeta_j^{\sigma}-\zeta_i^{\sigma}=\sigma(\delta_{ij})-\delta_{ij}\mbox{
  sur } U_i\cap U_j$.
La diff\'erentielle $\gd:C^{0,0}\rightarrow C^{1,0}\oplus C^{0,1}$ a pour
  image
$\gd\{\gamma_i\}=\{\{\sigma\gamma_i-\gamma_i\},\{\gamma_j-\gamma_i\}\}$.
Ainsi 
$$H^1(G,{\cal T}_C)={\{\{\zeta_i^{\sigma}\},\{\delta_{ij}\}\}\over \{\{\sigma\gamma_i-\gamma_i\},\{\gamma_j-\gamma_i\}\}}$$
Un 2-cocycle est un triplet $(\delta,\gamma,\beta)$ avec
$\{\delta_{ijk}\}$ un 2-cocycle de \v{C}ech, $\{\beta_i(\sigma,\tau)\}$
est un 2-cocycle au sens de la cohomologie des groupes (pour $i$
fix\'e) et $\{\gamma_{ij}(\sigma)\}$ est la composante de bidegr\'e
(1,1) ; donc
$$
^{\sigma}\delta_{ijk}-\delta_{ijk}=\gamma_{jk}(\sigma)-\gamma_{ik}(\sigma)+\gamma_{ij}(\sigma)=0$$
$$^\sigma\gamma_{ij}(\tau)-\gamma_{ij}(\sigma\tau)+\gamma_{ij}(\sigma)+\beta_j(\sigma,\tau)-\beta_i(\sigma,\tau)=0$$
\subsection{Identification de l'espace tangent et des obstructions}
Revenons \`a l'\'etude du foncteur des d\'eformations $D_{\glo}$ (voir
\S 2) et 
notons $t_{D_{\glo}}=D_{\glo}(k[\eps])$ l'espace tangent au foncteur
$D_{\glo}$.\\
\begin{prop}
\label{propidt}
Nous avons $t_{D_{\glo}}\cong H^1(G,{\cal T}_C)$.
\end{prop}
\begin{pr}
Il s'agit de d\'ecrire une d\'eformation de $(C,G)$ \`a $k[\eps]$, i.e. une d\'eformation $X$ de
$C$ \`a $k[\eps]$ et un rel\`evement de $G$ \`a $X$.
La d\'eformation $X$ est localement triviale, c'est-\`a-dire, il
existe un recouvrement ouvert affine $C=\cup_iU_i$, tel que
nous ayons des trivialisations locales de la d\'eformation $X$ de $C$
$$X_i\stackrel\sim\rightarrow U_i\otimes_k k[\eps]$$
non $G$-\'equivariantes. L'action de $G$ sur $X_i$ induit une action
sur $U_i\otimes_k k[\eps]$ qui est d\'ecrite par une classe de
cohomologie $\zeta_i=\{\zeta_{i,\sigma}\}_{\sigma\in G}\in
H^1(G,\Gamma(U_i,{\cal T}_C))$ ;
plus pr\'ecisement l'action de $G$ sur $\Gamma(U_i,{\cal T}_C)$ est l'action tordue par $\zeta_i :$
$$r\rightarrow r+\zeta_i(\sigma)(r)\eps$$
Sur $X_i\cap X_j$ il y a deux trivialisations locales de la structure de sch\'ema et les actions induites de $G$. 
Il y a compatibilit\'e de ces donn\'ees par l'isomorphisme $G$-\'equivariant de transition
$$\varphi_{ij}:(U_i\cap U_j)\otimes k[\eps]\rightarrow(U_i\cap U_j)\otimes k[\eps]$$
La $G$-\'equivariance s'\'ecrit
$\varphi_{ij}\zeta_j(\sigma)=\zeta_i(\sigma)^{\sigma}\varphi_{ij}$.
L'automorphisme $\varphi_{ij}$ est d\'ecrit par une d\'erivation 
$\{\delta_{ij}\}\in \Gamma(U_i\cap U_j,{\cal T}_C)$ :
$$\{\varphi_{ij}(r)\}=\{r+\delta_{ij}(r)\eps\}, \; r\in\Gamma(U_i\cap U_j,{\cal T}_C)$$
La condition de $G$-\'equivariance s'\'ecrit alors
$\delta_{ij}+\zeta_j(\sigma)=\zeta_i(\sigma)+ ^{\sigma}\delta_{ij}$ ;
soit encore
$$^{\sigma}\delta_{ij}-\delta_{ij}=\zeta_j(\sigma)-\zeta_i(\sigma)$$
Donc $\{\zeta_i(\sigma),\delta_{ij}\}$ d\'efinit une classe de cohomologie de $H^1(G,{\cal T}_C)$ vu la description pr\'ec\'edente de ce groupe.
Remarquons enfin que le changement de $\zeta_i(\sigma)$ et de
$\delta_{ij}$ par des 1-cobords, revient \`a changer les
trivialisations locales, donc donne
la m\^eme classe de cohomologie dans  $H^1(G,{\cal T}_C)$.
D'o\`u la proposition.
\end{pr}
Nous allons de la m\^eme fa\c{c}on identifier les obstructions au
rel\`evement de $(C,G)$ \`a des classes de 2-cohomologie.
Soit $A'\rightarrow A$ une petite extension de ${\cal C}$ et $(X,G)$
un repr\'esentant d'une d\'eformation $[(X,G)]\in D_{\glo}(A)$ \`a $A$ de $(C,G)$. Comme $H^2(C,{\cal T}_C)=0$ (vu que dim$_kC=1$), il n'y
a pas d'obtruction \`a relever le sch\'ema $X$ de $A$ \`a $A'.$ Il ne reste plus qu'\`a relever $G$, pour cela raisonnons localement.
Soit $\{U_i\}$ un recouvrement de $X$ par des ouverts affines $G$-stables.
Choisissons un rel\`evement de $A$ \`a $A'$ du sch\'ema affine $X_i$
(de support $U_i$) en $\widetilde{X}_i$ ; pour r\'esumer, nous avons la
situation suivante :
\begin{diagram}[labelstyle=\textstyle]
\widetilde{X}_i &\rTo &\mbox{Spec}A'\\
\uTo&&\uTo\\
X_i&\rTo & \mbox{Spec}A\\
\uTo&&\uTo\\
U_i&\rTo &\mbox{Spec}k\\
\end{diagram}
Au-dessus de $U_i\cap U_j$, nous avons deux d\'eformations de $X_i\cap X_j$ \`a $A'$ donc un isomorphisme (non unique)
$$\sigma_{ij}:\widetilde{X}_{j}|_{U_i\cap U_j}\stackrel\sim\longrightarrow\widetilde{X}_i|_{U_i\cap U_j}$$
qui se r\'eduit \`a l'identit\'e sur $A.$
L'automorphisme infinit\'esimal de $\widetilde{X}_i|_{U_i\cap U_j\cap U_k}$
$$\theta_{ijk}=\sigma_{ij}\sigma_{jk}\sigma_{ik}^{-1}\mbox{ sur } U_i\cap U_j\cap U_k$$
est d\'ecrit par une d\'erivation $\delta_{ijk}\in\Gamma(U_i\cap
U_j\cap U_k,{\cal T}_C)$ :
$$\theta_{ijk}=\Id+\eps \delta_{ijk}$$
Et $\{\delta_{ijk}\}$ est un 2-cocycle de \v{C}ech de classe nulle car $H^2(C,{\cal T}_C)=0$.
A pr\'esent il faut  relever l'action de $G$. Pour $\sigma\in G$, notons $\tilde{\sigma}_i$ un rel\`evement quelconque de $\sigma_i=\sigma|_{X_i}$.
Donc sur l'ouvert $X_i\cap X_j$ :
\begin{diagram}[labelstyle=\textstyle]
\widetilde{X}_{j}|_{X_i\cap X_j}&\rTo^{\sim}_{\sigma_{ij}}&\widetilde{X}_i|_{X_i\cap X_j}\\
\uTo^{\tilde{\sigma}_j}_{\sim}&&\uTo^{\tilde{\sigma}_i}_{\sim}\\
\widetilde{X}_j|_{X_i\cap
  X_j}&\rTo^{\sim}_{\sigma_{ij}}&\widetilde{X}_i|_{X_i\cap X_j}\\
\end{diagram}
Ce diagramme n'est a priori pas commutatif, introduisons un automorphisme infinit\'esimal de $\widetilde{X}_i|_{U_i\cap U_j}$ de ``d\'efaut'' $f_{ij}(\sigma)$ d\'efini par
$\sigma_{ij}\tilde{\sigma}_j=f_{ij}(\sigma)\tilde{\sigma}_i\sigma_{ij},\; \mbox{ sur } X_i\cap X_j$.
Pour $\sigma,\tau\in G$, la comparaison de $\tilde{\sigma}_i$ et $\tilde{\tau}_i$ avec $(\widetilde{\sigma\tau})_i$ introduit l'automorphisme infinit\'esimal de $\widetilde{X}_i$ de ``d\'eviation'' d\'efini par 
$\tilde{\sigma}_i\tilde{\tau}_i=g_i(\sigma,\tau)(\widetilde{\sigma\tau})_i$
et $\{g_i(\sigma,\tau)\}$ fournit un 2-cocycle de $G$ dans $\Gamma(U_i,{\cal T}_C)$.\\
Nous allons prouver que les donn\'ees locales
$\{\theta_{ijk}\},\{f_{ij}(\sigma)\},\{g_i(\sigma,\tau)\}$
d\'efinissent un 2-cocycle et donc une classe de cohomologie de $H^2(G,{\cal T}_C)$.
De plus si nous changeons les identifications locales, alors ce cocycle est
simplement alt\'er\'e par un cobord.
\begin{lem}
Nous avons la relation
$f_{ij}(\sigma)\; ^\sigma f_{ij}(\tau)g_i(\sigma,\tau)f_{ij}(\sigma\tau)^{-1}=\sigma_{ij}g_j(\sigma,\tau)\sigma_{ij}^{-1}$.
\end{lem}
\begin{pr}
En effet
$\sigma_{ij}\tilde{\sigma}_j\tilde{\tau}_j=\sigma_{ij}g_j(\sigma,\tau)(\widetilde{\sigma\tau})_j$
et
$$\sigma_{ij}\tilde{\sigma}_j\tilde{\tau}_j=f_{ij}(\sigma)\tilde{\sigma}_i\sigma_{ij}\tilde{\tau}_j=f_{ij}(\sigma)\tilde{\sigma}_if_{ij}(\tau)\tilde{\tau}_i\sigma_{ij}$$
$$=f_{ij}(\sigma) ^{\sigma}f_{ij}(\tau)\tilde{\sigma}_i\tilde{\tau}_i\sigma_{ij}=f_{ij}(\sigma) ^{\sigma}f_{ij}(\tau)g_i(\sigma,\tau)(\widetilde{\sigma\tau})_i\sigma_{ij}$$
$$=f_{ij}(\sigma) ^{\sigma}f_{ij}(\tau)g_i(\sigma,\tau)f_{ij}(\sigma\tau)^{-1}\sigma_{ij}(\widetilde{\sigma\tau})_j$$
D'o\`u la relation attendue sur $X_i\cap X_j$ :
$$\sigma_{ij}g_j(\sigma,\tau)=f_{ij}(\sigma) ^{\sigma}f_{ij}(\tau)g_i(\sigma,\tau)f_{ij}(\sigma\tau)^{-1}\sigma_{ij}$$
\end{pr}
Comme $g_j(\sigma,\tau)$ est un automorphisme infinit\'esimal et que le groupe de ces automorphismes est ab\'elien, isomorphe au groupe additif des champs de vecteurs, nous pouvons \'ecrire 
$\sigma_{ij}g_j(\sigma,\tau)\sigma_{ij}^{-1}=g_j(\sigma,\tau)$.
Nous obtenons ainsi la relation
$$f_{ij}(\sigma)\;^{\sigma}f_{ij}(\tau)f_{ij}(\sigma\tau)^{-1}=g_j(\sigma,\tau)g_i(\sigma,\tau)^{-1}$$
Le changement de notation $\gamma_{ij}=f_{ij}^{-1}$ permet de
reconna\^\i tre la relation de 2-cocyle pour le triplet
$\{\theta,f,g\}$.
Il reste entre les donn\'ees $\{\theta,f,g\}$ une autre relation qui
s'identifie avec la relation de 2-cobord : sur $X_i\cap X_j\cap X_k$,
nous avons les automorphismes :
$$\widetilde{X}_k\tilde{\sigma}_k\stackrel{\sigma_{jk}}\longrightarrow\widetilde{X}_j
\tilde{\sigma}_j\stackrel{\sigma_{ij}}\longrightarrow\widetilde{X}_i\tilde{\sigma}_i$$
Par cons\'equent
$\sigma_{ij}\sigma_{jk}\tilde{\sigma}_k=\theta_{ijk}\sigma_{ik}\tilde{\sigma}_k$.
Or
$$\sigma_{ij}\sigma_{jk}\tilde{\sigma}_k=\sigma_{ij}f_{jk}(\sigma)\tilde{\sigma}_j\sigma_{jk}=f_{jk}(\sigma)\sigma_{ij}\tilde{\sigma}_j\sigma_{jk}=f_{jk}(\sigma)f_{ij}(\sigma)\tilde{\sigma}_i\theta_{ijk}\sigma_{ik}$$
Et
$$\theta_{ijk}\sigma_{ik}\tilde{\sigma}_k=\theta_{ijk}f_{ik}(\sigma)\tilde{\sigma}_i\sigma_{ik}$$
D'o\`u la relation
$$f_{jk}(\sigma)f_{ij}(\sigma)f_{ik}(\sigma)^{1}\;^{\sigma}\theta_{ijk}\theta_{ijk}^{-1}=1$$
c'est-\`a-dire la relation de 2-cobord en bidegr\'e (1,2).
Enfin il est facile de voir que le changement des identifications
locales alt\`ere le cocycle par un 2-cobord, donc ne change pas sa
classe de cohomologie.
Nous avons ainsi \'etabli
\begin{prop} 
\label{propobs}
Soit $(X,G)$ un repr\'esentant d'un \'el\'ement de $D_{\glo}(A)$.
La donn\'ee d'un rel\`evement de la courbe $X$ de $A$ \`a $A'$ d\'efinit une
classe de cohomologie dans $H^2(G,{\cal T}_C)$ dite classe
  d'obstruction de la d\'eformation de $(X,G)$ de $A$ \`a $A'$.
Il existe une d\'eformation de $(X,G)$ de $A$ \`a $A'$ si et
  seulement si cette classe d'obstruction est nulle. 
\end{prop}
\subsection{Construction du morphisme local-global}
L'objectif de ce paragraphe est de relier le foncteur des
d\'eformations $D_{\glo}$ aux foncteurs des d\'eformations
infinit\'esimales d\'efinis au-dessus des points de ramification
sauvage (voir 1.1,1.2).\\

Soit $y_1,\ldots,y_r\in C/G$ les points images des points de
ramification sauvage de $C\rightarrow C/G$. Pour $1\leq i\leq r$, notons $x_i$ un point
au-dessus de $y_i$ et $G_{x_i}$ son stabilisateur.
Soit $D_{\loc}=\prod_i D_{G_{x_i}}$ o\`u $D_{G_{x_i}}$ est le foncteur des
d\'eformations infinit\'esimales de l'action induite de $G_{x_i}$ sur
$\widehat{\cal O}_{C,x_i}\cong k[[T]]$.
Le choix du point $x_i$ est sans importance, car si nous changeons
$x_i$ en autre point $x_i'$ de son orbite sous $G$, nous
obtenons un foncteur $D_{G_{x_i'}}$ isomorphe \`a $D_{G_{x_i}}$.
Par localisation en ces points, nous d\'efinissons un morphisme
$$\phi:D_{\glo}\rightarrow D_{\loc}$$
qui \`a une d\'eformation $X$ de $(C,G)$ \`a $A$ objet de ${\cal
  C}$ associe les d\'eformations infinit\'esimales  des
$G_{x_i}\hookrightarrow \Aut k[[T]]$ \`a $A$. Notons
$t_{D_{\loc}}=D_{\loc}(k[\eps])$ l'espace tangent au foncteur $D_{\loc}$.
\begin{lem}
\label{lem33311}
Le faisceau $R^1\pi_*^G({\cal T}_C)$ sur $C/G$ est concentr\'e aux points
$y_i$, et nous avons
$$R^1\pi_*^G({\cal T_C})_{y_i}\cong H^1(G_i,\widehat{{\cal T}}_{C,x_i})$$
\end{lem}
\begin{pr}
Le foncteur $R^n\pi_*^G({\cal T}_C)$ est ``local'', c'est-\`a-dire pour $V$ ouvert de
$C/G$
$$R^n\pi_*^G({\cal T}_C)|_V=R^n\pi^G_*|_{\pi^{-1}(V)}({\cal T}_C|_{\pi^{-1}(V)})$$
Nous nous restreignons alors \`a un ouvert affine $V$, c'est-\`a-dire
nous supposons que $C$ (et $C/G$) est
affine, ainsi
$$H^p(C/G,R^q\pi_*^G({\cal T}_C))=0, \mbox{ si } p>0,\mbox{ et }H^p(G,H^q(C,{\cal T}_C))=0, \mbox{ si } q>0$$
Donc les deux suites spectrales $'E_2^{p,q}$ et  $''E_2^{p,q}$
conduisent \`a :
$$H^0(C/G,R^q\pi_*^G({\cal T}_C))\cong H^q(G,{\cal T}_C)\cong
H^q(G,\Gamma(C,{\cal T}_C))$$
Donc dans le cas affine, les faisceaux $R^q\pi_*^G({\cal T}_C)$ et
$H^q(G,{\cal T}_C)\widetilde{\;\;}=H^q(G,\Gamma(C,{\cal
  T}_C))\widetilde{\;\;}$ co\"\i ncident.
Le foncteur de localisation \'etant exact, nous obtenons par passage \`a
la fibre en $y\in C/G$
$$R^q\pi_*^G({\cal T}_C)_y=H^q(G,\pi_*^G({\cal T}_C)_y)$$
Enfin comme $H^q(G,\pi_*^G({\cal T}_C)_y)$ est un module de longueur finie
(annul\'e par l'ordre de $G$), nous obtenons pour $x$ au-dessus de $y$, 
$$H^q(G,\pi_*^G({\cal T}_C)_y)\cong H^q(G,\pi_*^G({\cal T}_C)_y\widehat{\;\;}\;)
\cong H^q(G,\prod_{x\rightarrow y}\widehat{{\cal T}}_{C,x})\cong
H^q(G_x,\widehat{{\cal T}}_{C,x})$$
\end{pr}
\begin{lem}
L'application 
$H^1(G,{\cal T}_C)\rightarrow H^0(C/G, R^1\pi_*^G({\cal T}_C))\cong \oplus_{i=1}^r H^1(G_i,\hat{{\cal T}}_{C,x_i})$
est surjective  de noyau l'espace tangent au foncteur des
d\'eformations localement triviales.
\end{lem}
\begin{pr}
Consid\'erons la suite spectrale en degr\'e total $p+q=1$
$$'E_2^{p,q}=H^p(C/G,R^q\pi_*^G({\cal T}_C))\Rightarrow H^{p+q}(G,{\cal T}_C)$$
 Comme $'E_2^{2,0}=H^2(C/G,\pi_*^G({\cal
  T}_C))=0$, nous obtenons $'E_2^{0,1}='E_{\infty}^{0,1}$ ; 
nous avons \'egalement $'E_2^{1,0}='E_{\infty}^{1,0}$ ; d'o\`u la suite exacte
$$0\rightarrow H^1(C/G,\pi_*^G({\cal T}_C))\rightarrow H^1(G,{\cal T}_C)\rightarrow H^0(C/G,R^1\pi_*^G({\cal T}_C))\rightarrow 0$$
D'apr\`es le lemme \ref{lem33311}, 
$$H^0(C/G,R^1\pi_*^G({\cal T}_C))=\oplus_{i=1}^rH^1(G_i,\widehat{\cal T}_{C,x_i})$$
D'o\`u l'application surjective
$H^1(G,{\cal T}_C)\rightarrow \oplus_{i=1}^rH^1(G_i,\widehat{\cal T}_{C,x_i})$
qui est la diff\'erentielle de $\phi$, de noyau l'espace tangent au foncteur ker$\phi$ qui peut \^etre identifi\'e au foncteur des d\'eformations localement triviales.
\end{pr}
Passons \`a une analyse analogue du terme $H^2(G,{\cal T}_C)$.
\begin{lem}
Nous avons
$H^2(G,{\cal T}_C)\cong \oplus_{i=0}^r H^2(G_{x_i},\widehat{{\cal T}}_{C,x_i})$.
\end{lem}
\begin{pr}
Si $y\in C/G$ et $x\in\pi^{-1}(y)$, alors
$$R^2\pi_*^G({\cal T}_C)_y\widehat{\;\;}\cong H^2(G_x,\widehat{{\cal T}}_{C,x})$$
En effet cette fibre est nulle en dehors des points $y$ images de
points de ramification sauvage, donc le $({\cal T}_{C/G})_y$-module 
de longueur finie $R^2\pi_*^G({\cal T}_C)_y$ v\'erifie 
l'isomorphisme annonc\'e.
Dans la suite spectrale $'E_2^{p,q}=H^p(C/G,R^q\pi_*^G({\cal T}_C))$, nous avons
$'E_2^{2,0}=0$ car $H^2(C/G,\pi_*^G({\cal T}_C))=0$,
et $'E_2^{1,1}=0$ car $R^1\pi_*^G({\cal T}_C)$ est \`a support fini, donc
 $H^1(C/G,R^1\pi_*^G({\cal T}_C))=0$.
Ainsi par d\'eg\'en\'erescence, nous obtenons
$$'E_2^{0,2}=H^0(C/G,R^2\pi_*^G({\cal T}_C))=H^0(C/G,\oplus
R^2\pi_*^G({\cal T}_C)_y)\cong H^2(G,{\cal T}_C)$$
D'o\`u le r\'esultat annonc\'e.
\end{pr}
\begin{rem} Nous pouvons rapprocher ces r\'esultats de la proposition 1.5
  de Deligne Mumford \cite{DeMu} qui \'etablit un morphisme
  local-global similaire
pour le foncteur des d\'eformations de courbes stables, les points de
  branchement sauvage se substituant aux points doubles.
\end{rem}
\begin{thm}
\label{thlisse}
Le morphisme $\phi:D_{\glo}\rightarrow D_{\loc}$ est lisse.
\end{thm}
\begin{pr}
Soit $u:A'\rightarrow A$ une petite surjection de ${\cal C}$ ; nous avons le
diagramme suivant
\begin{diagram}[labelstyle=\textstyle]
D_{\glo}(A')&\rTo^{\phi} &D_{\loc}(A')\\
\uTo^{u}&&\uTo^{u}\\
D_{\glo}(A)&\rTo^{\phi} &D_{\loc}(A)\\
\end{diagram}
Soit $[(X,G)]\in D_{\glo}(A)$ tel que $\phi(X)$ se rel\`eve \`a $A'$ en $\widetilde{\phi(X)}$. Les obstructions \`a relever $X$ et $\phi(X)$ \`a
$A'$ co\"\i ncident, comme cette derni\`ere est nulle, nous avons donc une d\'eformation $\widetilde{X}$ de 
$(X,G)$ \`a $A'$.
Les \'el\'ements $\phi(\widetilde{X})$ et $\widetilde{\phi(X)}$ diff\`erent de
l'action d'un \'el\'ement de l'espace tangent $\delta\in t_{D_{\loc}}$; or $\phi$ est surjective sur les espaces tangents  
donc $\delta$ admet un ant\'ec\'edent par $\phi$ dans $t_{D_{\glo}}$.
Enfin si nous corrigeons $\widetilde{X}$ par cet ant\'ec\'edent nous obtenons l'\'egalit\'e
$\phi(\widetilde{X})=\widetilde{\phi(X)}$ ; d'o\`u la lissit\'e de $\phi$.
\end{pr}
\begin{thm}
\label{thmloglo}
Soit $R_i$ l'anneau de d\'eformations versel de $D_{G_{x_i}}$ et $R_{\glo}$
l'anneau de d\'eformations (uni)versel de $D_{\glo}$. Alors
$R_1\hat{\otimes}\cdots\hat{\otimes}R_r$ est l'anneau de d\'eformations versel
de $D_{\loc}$ et 
$$R_{\glo}=(R_1\hat{\otimes}\cdots\hat{\otimes}R_r)[[U_1,\ldots,U_N]]$$
avec $N=\dim_k H^1(C/G,\pi_*^G({\cal T}_C))$.
\end{thm}
\section{Le cas $G$ cyclique}
\label{gpcyc}
Nous supposons dor\'enavant que $G$ est un $p$-groupe cyclique.
L'objet de ce paragraphe est l'\'etude de l'anneau de d\'eformations
versel
local des d\'eformations d'une repr\'esentation $G\hookrightarrow\Aut
k[[T]]$. Nous commen\c{c}ons par d\'ecrire l'espace tangent et les
classes d'obstruction de ce probl\`eme de d\'eformations. Puis lorsque
$G$ est cyclique d'ordre $p$, nous mettons en \'evidence pour tout
$p>2$ une classe d'obstruction non triviale universelle.
\subsection{L'espace tangent}
Supposons \`a pr\'esent que $G$ est un groupe cyclique d'ordre
$p^n$.
 Fixons un g\'en\'erateur $\sigma$ de $G$ et notons 
$$\sigma\cdot T=f_{\sigma}(T)=T+\sum_{j\geq e}a_jT^j$$
avec $e>1$ et $a_e\not=0$ (l'action de $G$ est suppos\'ee non
triviale). Le groupe $G$ s'identifie alors au groupe de Galois de
l'extension $k((T))/k((Y))$ avec $k[[Y]]=k[[T]]^G$ la $k$-alg\`ebre des
invariants de $G$. Soit le $k[[T]]$-module libre $\displaystyle\Theta=k[[T]]{d\over
  dT}$ des $k$-d\'erivations continues. 
Il y a  une action naturelle du groupe $\Aut k[[T]]$ sur $\Theta$
et donc  une action de $G$ sur $\Theta$.
Les d\'eformations  de la repr\'esentation $G\hookrightarrow \Aut
k[[T]]$ sont donc d\'ecrites par les groupes de cohomologies
$H^1(G,\Theta)$, $H^2(G,\Theta)$. 
Nous voulons trouver dim$_kH^1(G,\Theta)$ et dim$_kH^2(G,\Theta)$ en
fonction des invariants de l'extension totalement ramifi\'ee
$k((T))/k((Y))$, c'est-\`a-dire de la filtration de $G$ par les
groupes de ramification sup\'erieurs (\cite{Se1} IV). 
Comme $G$ est cyclique, nous avons ([Se1] VIII 4) 
$$H^1(G,\Theta)={\Ker N\over\mbox{Im}\delta}\mbox{ et
  }H^2(G,\Theta)={\Ker\delta\over\mbox{Im}N}$$
pour $N,\delta$ op\'erateurs de $\Theta$, $k[[Y]]$-lin\'eaires d\'efinis par
$$\displaystyle N(h)=\sum_{i=0}^{p^n-1}\;^{\sigma^i}h,\;\; \delta(h)=\;^{\sigma}h-h,\;\;
h\in\Theta$$
Il s'agit d'abord d'\'etudier l'op\'erateur nilpotent d'ordre $p^n$,
$\delta=\sigma-1$, vu comme op\'erateur $k[[Y]]$-lin\'eaire de
$k[[T]]$ ou plus g\'en\'eralement de $E=T^{\beta}k[[T]]$ pour un $\beta\in
\N$ fix\'e.
Notons 
$$E_j=\Ker\delta^j=\{x \in E,\; \delta^jx =0\},\; 0\leq j\leq p^n$$
une suite d\'ecroissante de $k[[Y]]$ sous-modules sans torsion de
$E$. Comme $\delta^{-1}(E_j)=E_{j+1}$,
 $\delta$ induit une application $k[[Y]]$-lin\'eaire injective
$$\bar{\delta}: {E_{j+1}\over E_{j}}\longrightarrow {E_{j}\over E_{j-1}}$$
dont l'image est un $k[[Y]]$-module libre sans torsion.
Nous rappellons que $\lceil x\rceil$ d\'esigne la partie enti\`ere
sup\'erieure de $x$ et $\lfloor x \rfloor$ d\'esigne la partie
inf\'erieure de $x$.
\begin{thm} 
\label{thmbeta}
Soit $\beta$ l'exposant de la diff\'erente de
  $k((T))/k((Y))$. Alors
$$\dim_k H^1(G,\Theta)=\dim_kH^2(G,\Theta)=\Big\lfloor{2\beta\over p^n}\Big\rfloor-\Big\lceil{\beta\over p^n}\Big\rceil$$
\end{thm}
\begin{pr} Remarquons d'abord que
$$H^1(G,E)=\tors\Big({E\over\delta(E)}\Big)={E_{p^n-1}\over \delta(E)}$$
o\`u tors d\'esigne le sous-module de torsion.\\
En effet $\ker N$ est un sous-$k[[Y]]$-module de $E$ et le
$k[[Y]]$-module $E/\ker N$ est sans torsion.
La derni\`ere \'egalit\'e s'obtient en remarquant que 
$\delta(E)\subset E_{p^n-1}$ et si $x\in E$ v\'erifie $\delta y=lx$, alors
$\delta^{p^n-1}(lx)=l\delta^{p^n-1}(x)=0$ donc $x\in E_{p^n-1}$. Enfin
$E_{p^n-1}/\delta(E)$ est de torsion car $E_{p^n-1}$ et $\delta(E)$
sont des $k[[Y]]$-modules de m\^eme rang $p^n-1$.\\

Nous cherchons donc la dimension du $k$-espace vectoriel 
$\displaystyle V={E_{p^n-1}\over \delta(E)}$.
Ce $k[[Y]]$-module $V$ est filtr\'e par les sous modules
$$V^j={E_j+\delta(E)\over\delta(E)}\cong {E_j\over\delta(E_{j+1})},\; 1\leq j\leq p^n-1$$
et $V^0=0$ ; ainsi l'ordre de $V$ v\'erifie
$$|V|=\sum_{j=0}^{p^n-2}|V^{j+1}/V^j|=\sum_{j=0}^{p^n-2}|Q_j|$$
pour
$$Q_j\cong {E_{j+1}\over\delta(E_{j+2})+E_j}\cong\mbox{Coker}\Big(\bar{\delta}:{E_{j+2}\over E_{j+1}}\rightarrow{E_{j+1}\over E_j}\Big)$$
Donc nous avons
$$|V|=\Big|\mbox{Coker}\Big(E\cong{E_{p^n}\over
  E_{p^n-1}}\rightarrow {E_{p^n-1}\over E_{p^n-2}}\rightarrow
  \cdots\rightarrow {E_2\over E_1}\rightarrow E_1\Big)\Big|$$
d'o\`u $|V|=|\mbox{Coker}(E\rightarrow E_1)|$.
Nous avons l'\'egalit\'e entre op\'erateurs
$\delta^{p^n-1}=N=1+\sigma+\cdots+\sigma^{p^n-1}$.
Par cons\'equent
$$\displaystyle |V|={\mbox{dim}_kE^1\over\mbox{dim}_k\delta^{p^n-1}(E)}=
\mbox{dim}_k{\Ker \delta\over N(E)}$$
et ainsi 
$\dim_kH^1(G,E)=\dim_kH^2(G,E)$.
Pour obtenir les dimensions des groupes de
cohomologie $H^1(G,E)$, $H^2(G,E)$, il suffit donc d'\'etudier
$H^2(G,E)$.\\
Pour le calcul de $\dim_k H^2(G,\Theta)$, il est commode
d'identifier le $G$-module $\Theta$ \`a un sous-$G$-module de
$k[[T]]$. Soit ${\cal D}$ la diff\'erente de $k((T))/k((Y))$, 
$${\cal D}\cong (T^{\beta}) \mbox{ avec } \beta=\sum_{j=1}^\infty (|G_j|-1)$$
les $G_j$ \'etant les groupes de ramification sup\'erieurs
(\cite{Se1}).
Nous avons un isomorphisme de $G$-modules : 
$\Theta\cong {\cal D}$.
En effet,
posons $\displaystyle j_{\sigma}={df_{\sigma}\over dT}$ ; nous v\'erifions alors que
$\{j_{\sigma}\}$ est un cocycle de $G$ dans $k((T))^*$. D'apr\`es le
th\'eor\`eme 90 d'Hilbert c'est un cobord, c'est-\`a-dire   il existe $j\in k((T))^*$
tel que $\forall \sigma\in G,\;j_{\sigma}=j/\sigma j$. Nous pouvons prendre
$\displaystyle j={dY\over dT}$, $Y$ \'etant la norme de $T$. Par
d\'efinition $v_T(j)=\beta$.
Donc la multiplication par $T^{\beta}$ dans $k((T))$ \'etablit un
$G$-isomorphisme de $\Theta\cong k[[T]]$ avec le $G$-module
$T^{\beta}k[[T]]$ muni de l'action standard.\\ 
Nous avons 
$$H^2(G,E)={\Ker \delta\over N(E)}={E\cap k[[Y]]\over N(E)}$$
Or d'une part $E\cap k[[Y]]$ est le plus petit id\'eal de $k[[T]]$ contenant
$Y$ et $T^{\beta}$. Donc 
$$E\cap k[[Y]]=Y^{\lceil\beta/ p^n\rceil}k[[T]]$$
D'autre part nous connaissons l'image de la trace ([Se2] V 3) :\\
$$N(E)=Y^{\lfloor 2\beta/ p^n\rfloor}k[[T]]$$
D'o\`u le th\'eor\`eme annonc\'e.
\end{pr}
L'\'egalit\'e $\dim_k H^1(G,\Theta)=\dim_k H^2(G,\Theta)$ peut
\'egalement \^etre d\'eduite du th\'eor\`eme de Herbrand. Dans le
paragraphe suivant nous pr\'ecisons la structure de $H^1(G,\Theta)$
comme module sur $k[[Y]]=k[[T]]^G$.
\subsection{Une obstruction universelle}
Soit $R$ un anneau complet de valuation discr\`ete de corps
r\'esiduel $k$ et $k\subset R$ ; donc  
$R=k[[Y]]$. Soit
$K=k((Y))$ son corps des fractions. Consid\'erons l'\'equation
\begin{equation}
\label{equat}
X^p-X=\phi\in k((Y))\end{equation}
Il est bien connu que si $v_Y(\phi)=-m<0$, avec $p\wedge m=1$, alors
l'\'equation (\ref{equat}) est irr\'eductible et d\'efinie une
extension $p$-cyclique de conducteur de
Hasse $m$. Cela signifie que si $\{G_i\}$ est la filtration de $G$ par
les groupes de ramification sup\'erieurs, nous avons
$$G=G_0=\cdots=G_m\;\;\;\mbox{ et }\;\;\; G_{m+1}=\{1\}$$
Supposons dor\'enavant cette condition r\'ealis\'ee, et soit
$L=K(\xi)$, o\`u $\xi$ est une racine de (\ref{equat}). L'extension
$L/K$ est d'Artin-Schreier \index{extension ! d'Artin-Schreier} et $\eta=\xi$ est un g\'en\'erateur
d'Artin-Schreier,\index{g\'en\'erateur d'Artin-Schreier} c'est-\`a-dire $L=K(\eta)$ et $\eta^p-\eta\in
K$. Notons que tout g\'en\'erateur d'Artin-Schreier est de la forme
$\eta=j\xi+\psi$, o\`u $j\in\F_p^*$ et $\psi\in K$. Nous pouvons
parler plus g\'en\'eralement des g\'en\'erateurs d'Artin-Schreier
$\eta$ de conducteur donn\'e $m$ avec $\eta\in\overline{k((Y))}$, la
cl\^oture alg\'ebrique de $k((Y))$. Soit $\varphi\in \Aut k[[Y]]$ un automorphisme continu. Nous pouvons prolonger $\varphi$ \`a $\Aut k((Y))$, donc
$\varphi$ agit sur les g\'en\'erateurs d'Artin-Schreier. Nous avons le lemme :
\begin{lem}
Soit $\xi$ et $\eta$ deux \'el\'ements d'Artin-Schreier de conducteur
$m$. Alors il existe un automorphisme continu $\varphi\in\Aut k[[Y]]$
tel que $\varphi(\xi)=\eta$.
\end{lem}
\begin{pr}
Soit $\phi\in k((Y))$ tel que $v_Y(\phi)=-m$, avec $p\wedge m=1$. Commen\c{c}ons par d\'emontrer qu'il existe $Z\in k[[Y]]$ $v_Y(Z)=1$
(c'est-\`a-dire $Z$ est une variable) tel que $\phi=Z^{-m}$.
Il existe $\psi\in k[[Y]]^*$ avec $\phi=Y^{-m}\psi$. Comme
$m$ est premier \`a  $p$, $Y^m-\psi=0$ est une equation s\'eparable modulo $Y$.
Ainsi d'apr\`es le lemme d'Hensel, $T^m-\psi$ a $m$ racines distinctes. Soit
$P(Y)$ l'une de ces racines; alors $\phi=(Y^{-1}P(Y))^m=Z^{-m}$ avec
$Z=P(Y)^{-1}Y\in k[[Y]]$ une variable.
Donc il existe deux  variables $Z,W\in k[[Y]]$ telles que
$$\xi^p-\xi=Z^{-m}, \; \eta^p-\eta=W^{-m}$$
Soit $\varphi\in\Aut k[[Y]]$ tel que $\varphi(Z)=W$. Alors
$\varphi(\xi)$ est solution de $X^p-X=W^{-m}$ et il existe
$j\in\F_p$ avec $\varphi(\xi)=\eta+j$.
Soit $\sigma_{\eta}\in \Aut k[[Y]]$ avec $\sigma_{\eta}(\eta)=\eta+1$.
Nous avons $\sigma_{\eta}^{-j}\circ \varphi\in \Aut k[[Y]]$ et
$\sigma_{\eta}^{-j}\circ\varphi(\xi)=\eta$.
\end{pr}
Ainsi les automorphismes d'ordre $p$ et de conducteur $m$ de $k[[T]]$
sont deux \`a deux conjugu\'es et en particulier conjugu\'es \`a
$\sigma_0$, o\`u $\sigma_0(T)=T(1+T^m)^{-1/m}$. La proposition
suivante est analogue \`a 1.3(a) \cite{Ma} :
 \begin{prop}
\label{Maz}
Soit $\bar{\rho}$ et $\bar{\rho}'$ deux repr\'esentations de $G$ dans
$\Aut k[[T]]$ conjugu\'ees par un \'el\'ement de $\bar{\delta}\in\Aut
k[[T]]$. Alors les anneaux de d\'eformations versels
$R_{\bar{\rho}}$ et $R_{\bar{\rho}'}$ de $\bar{\rho}$ et $\bar{\rho}'$ respectivement sont isomorphes.
\end{prop}
Pour \'etudier les anneaux de d\'eformations versels de repr\'esentations
d'un groupe cyclique d'ordre $p$ dans $\Aut k[[T]]$, il suffit donc
d'\'etudier les d\'eformations de l'automorphisme d'ordre $p$
$$\sigma_0(T)={T\over (1+T^m)^{1/m}}, \;\;\; p\wedge m=1$$
\begin{lem}
\label{eltH1}
L'\'el\'ement
$$h(T){d\over dT}=\left\{\begin{array}{cc}{d\over dT}&\mbox{ si }
    m=1\cr
T{d\over dT}&\mbox{ si } m>1\end{array}\right.$$
repr\'esente un \'el\'ement non nul de $H^1(G,\Theta)$.
\end{lem}
\begin{pr}
Rappelons que $\displaystyle H^1(G,\Theta)\cong {\Ker N\over \im
  \delta}$ et 
$\displaystyle N(h(T))=\sum_{i=0}^{p-1}\sigma_0^i(h(T))\Big({d\sigma_0^i(T)\over
  dT}\Big)^{-1}$. 
Comme $\displaystyle \Big({d\sigma_0^i(T)\over
  dT}\Big)^{-1}=(1+iT^m)^{1+1/m}$, nous avons si $m>1$
$$N(h(T))=N(T)=\sum_{i=0}^{p-1}{T\over (1+iT^m)^{1/m}}(1+iT^m)^{1+1/m}=\sum_{i=0}^{p-1}T+iT^{m+1}=0$$
et si $m=1$
$$N(h(T))=N(1)=\sum_{i=0}^{p-1}(1+iT)^{2}=\sum_{i=0}^{p-1}1+2iT+i^2T^2=0$$
Donc $h(T)\in \ker N$.
Pour $f\in k[[T]]$, $v_T(\delta(f))\geq m$, donc $h(T)\not\in\im\delta.$
\end{pr}
Soit $A$ un objet de ${\cal C}$.
D\'eformer suivant la direction tangentielle d\'efinie par
$\displaystyle h(T){d\over dT}$, conduit \`a \'etudier les
automorphismes suivants :
$$\sigma_a(T)={T\over(1+aT^m)^{1/m}} \mbox{ pour } a\in 1+{\cal M}_A \mbox
{ si } m>1 $$
$$\sigma_a(T)={T+a\over 1+T+a} \mbox{ pour } a\in {\cal M}_A \mbox
{ si } m=1 $$
Commen\c{c}ons par le cas $m>1$, $p>2$,  et l'observation \'el\'ementaire
\begin{lem}
\label{defsig}
Soit $A$ objet de ${\cal C}$, $a\in 1+{\cal M}_A$ et 
$$\sigma_a(T)=T/(1+aT^m)^{1/m}$$
Nous avons $\sigma_a^p=\Id$ si et seulement si $1+a+\cdots+a^{p-1}=0$.
\end{lem}
\begin{lem}
\label{condobs}
Supposons $p>2$, $m>1$ tels que $p\wedge m=1$.
Soit $A'\rightarrow A$ une petite extension et $a\in{\cal M}_A$ tel
que $\sigma_a^p=\Id$. La classe d'obstruction  \`a relever $\sigma_a$ de $A$ \`a $A'$ est nulle
si et seulement si $1+a'+\cdots+{a'}^{p-1}=0$.
\end{lem}
\begin{pr}
Il s'agit de relier la condition du lemme \ref{defsig} \`a une obstruction
cohomologique, c'est-\`a-dire \`a une classe de $H^2(G,\Theta)$.
Soit $\pi : A'\rightarrow A$  une petite extention de noyau $tA'$, et
$a\in 1+{\cal M}_A$ tel que $\sigma_a^p=\Id$. Pour obtenir
l'obstruction cohomologique \`a relever $\sigma_a$ de $A$ \`a $A'$, nous
relevons $a$ en $a'\in A$. Ainsi
$$1+a'+\cdots+{a'}^{p-1}\in tA'\mbox{ et } {a'}^p=1$$
et l'obstruction cohomologique dans $H^2(G,\Theta)$ est d\'etermin\'ee
par le cocycle
$(\varphi_{i,j})_{i,j\in\F_p}$ d\'efini par
$$\sigma_{a'}^i\sigma_{a'}^j=\varphi_{i,j}\sigma_{a'}^{i+j}$$
les exposants \'etant pris modulo $p$. Donc
$$\varphi_{i,j}=\left\{\begin{array}{ll}1&\mbox{ si } i+j<p\\
\sigma_{a'}^p&\mbox{ si } i+j\geq p
\end{array}\right. \;\;\;\;i,j<p$$
L'exposant de la diff\'erente est $\beta=(m+1)(p-1)$ et
$H^2(G,\Theta)$ est le $k[[Y]]$-module cyclique
$$H^2(G,\Theta)={\Ker\delta\over
  \im N}\cong{(Y)^{\lceil {\beta/ p}\rceil}\over (Y)^{\lfloor
  2\beta/p\rfloor}}$$
$Y$ \'etant la norme de $T$ et l'isomorphisme de la preuve du
th\'eor\`eme \ref{thmbeta}
$$h(T){d\over dT}\mapsto h(T){dY\over dT}$$
\'etant \'etendu aux classes de cohomologie. Dans cette identification,
$\displaystyle \varphi_{i,j}=h_{i,j}{d\over dT}$ correspond \`a 
$$h_{i,j}=\left\{\begin{array}{ll}0&\mbox{ si } i+j<p\\
\sigma_{a'}^p(T)-T&\mbox{ si } i+j\geq p
\end{array}\right. \;\;\;\;i,j<p$$
D'o\`u la classe d'obstruction est la classe de
$\displaystyle
h_{p-1,1}{d\over dT}=(\sigma_{a'}^p(T)-T){d\over dT}$. Pour v\'erifier
  que cette classe est non nulle, observons que,
$$\sigma_{a'}^p(T)-T={T\over
  ({{a'}^p}+(1+a'+\cdots+{a'}^{p-1})T^m)^{1/m}}-T=
-{1\over m}(1+a'+\cdots+{a'}^{p-1})T^{m+1}$$
Il suffit alors de v\'erifier  
\begin{equation}
\label{equineg}
v_Y\Big(h_{p-1,1}(T){dY\over dT}\Big)<\Big\lfloor {2\beta\over
  p}\Big\rfloor
\end{equation}
Cette in\'egalit\'e est \'equivalente \`a
$$v_T(h_{p-1,1})+\beta=(m+1)p<p\lfloor2(m+1)(p-1)/p\rfloor$$
\'Ecrivons $2(m+1)=pu+u'$ avec $0\leq u'<p$.
Nous avons 
$$\lfloor 2(m+1)(p-1)/p \rfloor=\left\{\begin{array}{ll}u(p-1)&\mbox{
      si } u'=0\\
u(p-1)+u'-1 &\mbox{ sinon } 
\end{array}\right.$$
Si $u'\geq 1$, le r\'esultat est clair. Si $u'=0$ l'in\'egalit\'e est
assur\'ee pourvu que $p>2$ et  $(m,p)\not=(1,3)$.
Ce qui d\'emontre l'in\'egalit\'e (\ref{equineg}) et le lemme.
\end{pr}
Revenons \`a pr\'esent au cas $m=1$, $p>2$.
\begin{lem}
Soit $A$ objet de ${\cal C}$, $a\in {\cal M}_A$,
$$ \sigma_a(T)={T+a\over 1+T+a}\;\;\mbox{ et }\;\; M_a=\left(\begin{array}{cc}1&a\cr 1&1+a\end{array}\right)$$
Nous avons $\sigma_a^p=\Id$ si et seulement si $M_a^p=\Id$.
\end{lem}
\begin{pr}
Soit $A_p,B_p,C_p,D_p\in \Z[a]$ tels que
$$M_a^p=\begin{pmatrix}A_p&B_p\cr C_p&D_p\end{pmatrix}$$ 
L'automorphisme $\sigma_a$ s'identifie \`a l'homographie de matrice
$M_a$. Ainsi
$$\sigma_a^p=\Id\Longleftrightarrow \sigma_a^p(T)=T\Longleftrightarrow
A_pT+B_p=T(C_pT+D_p)$$
$$\sigma_a^p=\Id \Longleftrightarrow A_p=D_p,
B_p=C_p=0\Longleftrightarrow
M_a^p=\alpha\Id, \mbox{ pour }\alpha\in\Z[a]$$
Supposons $M_a^p=\alpha\Id$ pour $\alpha\in \Z[a]$. Comme $M_a \mbox
{ mod }{\cal M}_A=M_0$, la matrice $M_a^p$ rel\`eve $M_0^p=\Id$ \`a
$A$.
Donc $\alpha\equiv 1 \mbox{ mod }{\cal M}_A$.
Par ailleurs det$M_a=1$ donc det$M_a^p=\alpha^2=1$. Ces deux
conditions sur $\alpha$ entra\^\i nent que $\alpha=1$. D'o\`u
$$\alpha_a^p=\Id\Longleftrightarrow M_a^p=\Id$$ 
\end{pr}
Soit $T_p,S_{p-1}$ les polyn\^omes de Tchebychev
respectivement de premi\`ere et de seconde esp\`ece
$$T_p(X)={1\over 2}\sum_{l=0}^{\lfloor
  p/2\rfloor}\left(\begin{array}{c}p-l\cr l\end{array}\right){p\over
  p-l}(-1)^l(2X)^{p-2l}$$
$$S_{p-1}(X)=\sum_{l=0}^{
  (p-1)/2}\left(\begin{array}{c}p-1-l\cr
    l\end{array}\right)(-1)^l(2X)^{p-1-2l}$$
Par d\'efinition, si $2X=Z+Z^{-1}$, nous avons
$$2T_p(X)=Z^p+Z^{-p},\;\;\;\; (Z-Z^{-1})S_{p-1}(X)=Z^p-Z^{-p}$$
\begin{lem}
\label{lemzz}
Nous avons
$$M_a^p=\Id \Longleftrightarrow S_{p-1}(a/2+1)=T(a/2+1)-1=0$$
\end{lem}
\begin{pr}
Posons $\displaystyle z={2+a+\sqrt{(2+a)^2-4}\over 2}$. 
En diagonalisant $M_a$, nous
obtenons
$$\displaystyle M_a^p=\left(\begin{array}{cc}-{z^p-z^{-p}\over z-z^{-1}}{a\over 2}
    +{z^p+z^{-p}\over 2}& {z^p-z^{-p}\over
      z-z^{-1}}\Big(\Big({z-z^{-1}\over 2}\Big)^2-{a^2\over 4}\Big)\cr
{z^p-z^{-p}\over z-z^{-1}}& {z^p-z^{-p}\over z-z^{-1}}{a\over
  2}+{z^p+z^{-p}\over 2}\cr\end{array}\right)$$
Donc $M_a^p=\Id$ si et seulement si
$$ {z^p-z^{-p}\over z-z^{-1}}=0\;\;\mbox{ et }\;\;z^p+z^{-p}=2$$
A l'aide des polyn\^omes de Tchebychev de premi\`ere et seconde
esp\`ece ces conditions s'\'ecrivent
$$ S_{p-1}(x)=0 \mbox{ et }T_p(x)-1=0 \mbox{ pour } 2x=z+z^{-1}=2+a$$ 
\end{pr}
\begin{lem}
\label{lempsi}
L'id\'eal $(T_p(X/2+1),S_{p-1}(X/2+1))$ de $W(k)[X]$ est 
principal et engendr\'e par
$$\psi(X)=\sum_{l=0}^{
  (p-1)/2}\left(\begin{array}{c}p-1-l\cr l\end{array}\right)(-1)^l(X+4)^{(p-1)/2-l}$$
\end{lem}
\begin{pr}
En \'etudiant les racines de $T_p(X)-1$ et $S_{p-1}(X)$, vus comme polyn\^omes de $\C(X)$ nous pouvons
montrer
$$\displaystyle \phi(X)=(T_p(X)-1)\wedge S_{p-1}(X)=2^{(p-1)/2}\prod_{l=1}^{
  (p-1)/2}(X-\cos {2l\pi\over p})$$
Pr\'ecisement
$$\phi(X)=S_{p-1}\Big(\sqrt{{X+1\over 2}}\Big)=\sum_{l=0}^{
  (p-1)/2}\left(\begin{array}{c}p-1-l\cr l\end{array}\right)(-1)^l(2(X+1))^{(p-1)/2-l}$$
Enfin \`a l'aide de la trigonom\'etrie hyperbolique, nous obtenons une
\'equation de B\'ezout
$$U(X)(T_p(X)-1)+V(X)S_{p-1}(X)=\phi(X)$$
avec 
$\displaystyle U(X)=-{1\over 2}\phi(X)$ et,
$$V(X)={1\over 4}\sum_{l=0}^{\lfloor
  p/2\rfloor}\left(\begin{array}{c}p-l\cr l\end{array}\right){p\over
  p-l}(-1)^l(2(X+1))^{(p+1)/2-l}$$
Comme les entiers premiers \`a $p$ sont inversibles dans $W(k)$,
cette relation de B\'ezout est encore valable dans $W(k)[X]$, m\^eme
apr\`es le changement de variables $Y=X/2+1$.
Donc
l'id\'eal  $(T_p(X/2+1),S_{p-1}(X/2+1))$ de $W(k)[X]$ est 
principal et engendr\'e par
$\psi(X)=\phi(X/2+1)$. 
\end{pr}

Soit $A'\rightarrow A$ une petite extension de ${\cal C}$ de noyau
$tA'$, $a\in{\cal M}_A$ tel que $\sigma_a^p=\Id$ et $a'$ un
rel\`evement de $a$ dans $A'$. La classe d'obstruction \`a relever
$\sigma_a$
de $A$ \`a $A'$ est la classe de $\displaystyle h(T){d\over dT}$
dans $H^2(G,\Theta)$, avec
$$\sigma_{a'}^p(T)-T=-{a'T^2\over 1+a'T}=th(T)$$
Par r\'ecurrence sur $n\in \N$, nous pouvons montrer que $M_{a'}^n$ est de la
forme 
$$M_{a'}^n=\begin{pmatrix}A_n & a'C_n\cr C_n& A_n+a'C_n\end{pmatrix}$$
Ainsi
$$\sigma_{a'}^p(T)-T={A_pT+a'C_p\over C_pT+A_p+a'C_p}-T={a'C_p-a'C_pT-C_pT^2\over C_pT+A_p+a'C_p}$$
Or $C_p\equiv 0\mbox{ mod }{\cal M}_A$, donc $C_p=tc_p$ avec $c_p\in
k$ et $a'C_p=0$.\\
De m\^eme   $A_p\equiv 1\mbox{ mod }{\cal M}_A$,
donc $A_p=1+ta_p$ avec $a_p\in k$. D'o\`u
$$\sigma_{a'}^p(T)-T={-c_pT^2t\over 1+ta_p+tc_pT}=-c_pT^2(1-ta_p-tc_pT)=(-c_pT^2)t$$
Donc la classe d'obstruction est la classe de 
$\displaystyle -c_pT^2{dY\over dT}\in {(Y)^{\lceil \beta/ p\rceil}\over
  (Y)^{\lfloor 2\beta/p\rfloor}}$.
Par cons\'equent cette classe est nulle si et seulement si $c_p=0$.
Remarquons que si $c_p=0$, alors det$M_{a'}^p=1+2ta_p=1$ donc $a_p=0$. 
Donc nous avons montr\'e que l'obstruction cohomologique est nulle si
et seulement si $M_{a'}^p=\Id$. Pour r\'esumer nous avons le lemme
suivant (analogue du lemme \ref{condobs})
\begin{lem}
\label{condobs2}
Supposons $m=1$, $p>2$.
Soit $A'\rightarrow A$ une petite extension de ${\cal C}$, $a\in{\cal M}_A$ tel que $\sigma_a^p=\Id$ et $a'$ un
rel\`evement de $a$ \`a $A'$.
La classe d'obstruction \`a relever $\sigma_a$ de $A$ \`a $A'$ est
nulle si et seulement si $T_p(a'/2+1)-1=S_{p-1}(a'/2+1)=0$.
\end{lem}
Ces r\'esultats pr\'eliminaires \'etant \'etablis, nous pouvons d\'emontrer
\begin{thm}
\label{resultcond}
Soit $p>2$, $m\geq 1$ tels que $p\wedge m=1$. Soit
$$\sigma_0 : G\rightarrow \Aut k[[T]]$$
une repr\'esentation de conducteur de Hasse $m$.
Si $m>1$ alors nous avons un morphisme surjectif
$$R_{\sigma_0}\twoheadrightarrow W(k)[[X]]/\Big(\sum_{j=0}^{p-1}(1+X)^{mj}\Big)$$
Si $m=1$ et $p>3$, 
$R_{\sigma_0}\cong W(k)[[X]]/ (\psi(X))$. En particulier si $m=1$ et $p>3$ ou si $m=2$ et $p=5$, l'anneau de
d\'eformations versel $R_{\sigma_0}$ est d'intersection compl\`ete.
Enfin si $m=1$ et $p=3$ le probl\`eme de d\'eformations 
est rigide.
\end{thm}
\begin{pr} Si $(m,p)=(1,3)$, $\dim_kH^1(G,\Theta)=\dim_kH^2(G,\Theta)=0$, le
probl\`eme de d\'eformations est donc rigide. Pour $(m,p)\not=(1,3)$ ce r\'esultat s'obtient gr\^ace \`a
  l'\'etude du mod\`ele local $\sigma_0(T)=T(1+T^m)^{-1/m}$ puis la
  proposition \ref{Maz} conclut.\\
Commen\c cons l'\'etude du mod\`ele local par le cas $m=1$, $p>3$.\\
D'apr\`es les lemmes \ref{lemzz}, \ref{lempsi}, l'automorphisme
$$\sigma_X(T)={T+X\over 1+T+X}$$
d\'efinit une d\'eformation de $\sigma_0$ \`a $W(k)[[X]]/(\psi(X))$.  
Par d\'efinition de l'anneau versel, nous avons un morphisme
$$u:R_{\sigma_0}\rightarrow W(k)[[X]]/(\psi(X))$$
Comme 
$$\dim_k H^1(G,\Theta)=\lfloor 4(p-1)/p\rfloor -\lceil
2(p-1)/p\rceil=1$$
le morphisme $u$ induit un isomorphisme sur les espaces
tangents de $W(k)[[X]]/(\psi(X))$ et de $R_{\sigma_0}$. Nous allons montrer que $u$ est un isomorphisme en
d\'emontrant que le couple $(W(k)[[X]]/(\psi(X)),\sigma_X)$ d\'efinit une d\'eformation
verselle. 
Notons $D$ le foncteur de d\'eformations de $\sigma_0$ et
$h_{W(k)[[X]]/(\psi(X))}$ le foncteur des homomorphismes
$$h_{W(k)[[X]]/(\psi(X))}=\Hom_{W(k)}(W(k)[[X]]/(\psi(X)), \cdot)$$ 
Il suffit donc de v\'erifier que le
morphisme de foncteurs
$$D\rightarrow h_{W(k)[[X]]/(\psi(X))}$$
est lisse. Soit $A'\rightarrow A$ une petite
extension de ${\cal C}$ et $a\in {\cal M}_A$ tel que
$\sigma_a^p=\Id$. Si $[\sigma_a]\in D(A)$ se rel\`eve dans $D(A')$
alors, par d\'efinition, l'obstruction \`a relever $\sigma_a$ de $A$
\`a $A'$ est nulle. Alors pour $a'$ relevant $a$ \`a $A'$,
$\sigma_{a'}^p=\Id$ (Lemme \ref{condobs2}). Donc le morphisme de foncteurs
est lisse, et
$$R_{\sigma_0}\cong W(k)[[X]]/(\psi(X))$$
Un raisonnement essentiellement identique permet de traiter le cas $m>1$. 
\end{pr}
\begin{rem}
Si $m=1$ et $p>3$, $R_{\sigma_0}$ est d'intersection compl\`ete et
$R_{\sigma_0}/pR_{\sigma_0}=k[X]/(X^{(p-1)/2})$.
En effet d'apr\`es le lemme \ref{lempsi}, les racines de
$\psi(X)=\phi(X/2+1)$ sont les $\displaystyle 2(\cos{2l\pi\over p}-1)$
pour $l\in\{1,\ldots, (p-1)/2\}$. Donc $\psi(X)$, vu comme polyn\^ome
de $\Q[X]$, est le polyn\^ome
minimal de $e^{2i\pi/p}+ e^{-2i\pi/p}-2\in \Z[e^{2i\pi/p}]$. Or
 $\Q(e^{2i\pi/p}+ e^{-2i\pi/p})$ est une extension de $\Q$ totalement
 ramifi\'ee en $p$ d'indice $(p-1)/2$. Donc tous les coefficients de
 degr\'e inf\'erieur \`a $(p-1)/2$ de $\psi(X)$ sont divisible par
 $p$.
\end{rem}
\begin{ex}
Supposons $(p,m)=(5,2)$. Dans ce cas
$\dim_k H^1(G,\Theta)=1$. D'apr\`es le th\'eor\`eme  \ref{resultcond} l'anneau de d\'eformations
versel de l'automorphisme d'ordre $p$ d\'efini par
$\sigma_0(T)=T/(1+T^2)^{1/2}$ est
$$R_{\sigma_0}\cong  W(k)[[X]]/\Big(\sum_{j=0}^{p-1}(1+X)^{mj}\Big)$$
Donc $R_{\sigma_0}$ est d'intersection compl\`ete.
\end{ex}
 \subsection{D\'eformations des \'equations d'Artin-Schreier}
\label{subdefor}
  \'Ecrivons $m=pq-l$, $l\in[1,p[$.
  Soit $\zeta$ une racine de l'unit\'e d'ordre $p$. Posons
  $S=W(k)(\zeta)$, $\pi$ une uniformisante de $S$.
  Nous construisons une d\'eformation lisse d'un automorphisme
  $\sigma_0$ d'ordre $p$
  et de conducteur $m$ avec $n$ param\`etres ind\'ependants sur $S$, o\`u 
$$n=\left\{\begin{array}{ll}q&\mbox{  si } l=1\\
q-1&\mbox{ sinon }\end{array}\right.$$
Dans cette construction nous utilisons une id\'ee de
  Sekiguchi, Oort et Suwa \cite{SeOoSu}, exploit\'ee par Green et Matignon
  \cite{GrMa2}
  qui permet de voir l'\'equation d'Artin-Schreier (d\'efinissant
  l'automorphisme $\sigma_0$) comme la
  r\'eduction modulo $\pi$ d'une \'equation de Kummer. En d\'eformant
  ces \'equations, c'est-\`a-dire en prenant $A=S[[x_1,\ldots,x_n]]$
  pour anneau de coefficients, nous d\'efinissons une extension $L$ du
  corps des fractions de $A[[t]]$. Pour obtenir la
  d\'eformation recherch\'ee, il suffit d'identifier la 
normalisation de $A[[t]]$ dans $L$, qui est naturellement munie d'un
  automorphisme d'ordre $p$, \`a une alg\`ebre de s\'eries formelles. Cette
  strat\'egie est valable pour tout $l\in[1,p[$ cependant les cas $l=1$
  et $l\not=1$ doivent \^etre \'etudi\'es s\'epar\'ement.\\
  Soit $\lambda=\zeta-1$ ; donc
  $p\lambda^{-j}\equiv 0 \;(\mbox{mod} \pi),\mbox{ si } 0\leq j<p-1$,
  $p\lambda^{p-1}\equiv -1 \;(\mbox{mod} \pi)$.
  Consid\'erons l'\'equation de Kummer
  \begin{equation}
  \label{equaKum}
  (X+\lambda^{-1})^p-\lambda^{-p}=t^{-m}
  \end{equation}
  qui d\'efinit une extension cyclique du corps des fractions $K$ de $S[[t]]$.
  L'\'equation (\ref{equaKum}) est \'equivalente \`a
  \begin{equation}
  \label{equadev}
  X^{-p}=t^m\Big(1+{p\over \lambda}X^{-1}+\cdots +{p\over
  \lambda^{p-1}}X^{-(p-1)}\Big)
  \end{equation}
  Cette \'equation montre que $X^{-1}$ appartient \`a la normalisation
  $B$ de $S[[t]]$ dans $L=K(X)$. Par r\'eduction modulo $\pi$ l'\'equation
  (\ref{equadev}) conduit \`a l'\'equation d'Artin-Schreier de conducteur $m$
  \begin{equation}
  \label{Artin}
  X^{p}-X=t^{-m}
  \end{equation}
d\'efinissant l'automorphisme $\sigma_0$ de $k[[t]]$ par $\sigma_0(X)=X+1$.\\
 Commen\c{c}ons par \'etudier le cas $l=1$. Posons $\xi=Xt^q$. L'\'equation
  d'Artin-Schreier (\ref{Artin}) s'\'ecrit
  \begin{equation}
  \label{equal1}
  \xi^p-t^{(p-1)q}\xi=t
  \end{equation}
  et l'\'equation (\ref{equaKum}) 
  \begin{equation}
  \label{equaxi}
  (\lambda\xi+t^q)^p-t^{qp}=\lambda^pt
  \end{equation}
\begin{lem}
  \label{lemnorm}
  La normalisation $B$ de $S[[t]]$ dans $L$ est l'anneau des s\'eries formelles $S[[\xi]]$.
  \end{lem}
  \begin{pr}
  L'\'equation (\ref{equaxi}) \'equivaut \`a
  \begin{equation}
  \label{equapr}
  \xi^p+{p\over \lambda}\xi^{p-1}t^q+\cdots+{p\over\lambda^{p-1}}\xi
  t^{(p-1)q}=t
  \end{equation}
  Donc $\xi\in B$ et (\ref{equapr}) est de la forme $F(\xi,t)=0$ avec $
  {\partial F\over \partial t}=1\mod(\xi,t)$ ; nous pouvons alors r\'esoudre cette
  \'equation implicite en $t$ et plonger $S[[t]]$ dans $S[[\xi]]$, ce qui conduit
  \`a $B=S[[\xi]]$.
  \end{pr}
  Nous d\'eformons l'\'equation (\ref{equal1}) par la substitution de $t^q$
par $a(t)=t^q+x_1t^{q-1}+\cdots+x_q$. Cela signifie que nous consid\`erons \`a
pr\'esent l'anneau des s\'eries formelles $A=S[[x_1,\ldots,x_q]]$ et de la
m\^eme mani\`ere nous d\'eformons (\ref{equapr}) en 
\begin{equation}
\label{equadefo}
t=\xi^p+{p\over\lambda}\xi^{p-1}a(t)+\cdots+{p\over\lambda^{p-1}}\xi a(t)^{p-1}
\end{equation}
qui modulo $\pi$ se r\'eduit \`a
\begin{equation}
\label{equadered}
t=\xi^p-\xi a(t)^{p-1}
\end{equation}
Le lemme suivant se d\'emontre de mani\`ere analogue au lemme \ref{lemnorm}
\begin{lem}
\label{aririri}
Soit $L$ l'extension du corps des fractions de $A$ d\'efinie par
(\ref{equadefo}). Alors la normalisation de $A[[t]]$ dans $L$ est
$B=A[[\xi]]$.
\end{lem}
Pour traiter le cas $l\not=1$, notons le lemme suivant :
\begin{lem}
\label{lemBez}
Soit $B$ un anneau int\'egralement clos de corps des
fractions $L$. Soit $u$ et $v$ \'elements de $B$, $r,s\in \N$ avec
$r\wedge s=1$ et $u^r=v^s$. Alors il existe $\xi\in B$ tel que $u=\xi^s$. 
\end{lem}
\begin{pr}
Soit les coefficients de B\'ezout $cr+ds=1$. Nous avons
$u^{cr}=u^{1-ds}=v^{cs}$, donc $u=\xi^s$ avec $\xi=v^cu^d\in L$. Comme
$B$ est int\'egralement clos, $\xi\in B.$\end{pr}
Supposons d'abord $A=S$ ; il r\'esulte de l'\'equation (\ref{equadev})
que nous pouvons extraire une racine $m$-\`eme 
$Z=X^{-1/m}$ dans ${\cal M}_B$ et
\begin{equation}
\label{equaZ}
Z^p=t\Big(1+{p\over\lambda}X^{-1}+\cdots+{p\over\lambda^{p-1}}X^{-(p-1)}\Big)^{1/m}
\end{equation}
Si $L$ est l'extension du corps des fractions $K$ de $A[[t]]$
d\'efinie par l'\'equation (\ref{equadev}), nous d\'eduisons de
(\ref{equaZ}) que $t\in A[[Z]]$ et donc que la
normalisation $B$ de $A[[t]]$ dans $L$ est $A[[Z]]$ 
(comparer au th\'eor\`eme
II 4.1 de \cite{GrMa1} qui utilise un crit\`ere local de bonne
r\'eduction pour un rel\`evement de l'\'equation d'Artin-Schreier 
\`a un anneau de valuation discr\`ete).
Notons que par r\'eduction modulo $\pi$, l'inclusion $S[[t]]\subset
S[[Z]]$ donne la normalisation $k[[Z]]$ de $k[[t]]$ dans l'extension
d'Artin-Schreier $k((t))(X)$, $X^p-X=t^{-m}$ d\'efinie par l'\'equation
(\ref{Artin}), avec $Z=X^{-1/m}$.\\
Pour d\'eformer ce rev\^etement, nous pouvons r\'ep\'eter tel quel
l'argument pr\'ec\'edent. Comme pour le cas $l=1$, changeons le
param\`etre $X$ en $\xi=Xt^q$.
Supposons maintenant que l'anneau de base est
$A=S[[x_1,\ldots,x_{q-1}]]$. Soit
$a(t)=t^q+x_1t^{q-1}+\cdots+x_{q-1}t$. Par rapport au cas $l=1$, nous
annulons le terme constant. En effet contrairement au cas $l=1$
l'\'equation (\ref{equadev}) d\'eform\'ee  
\begin{equation}
\label{equadefol}
t^l=\xi^p+{p\over\lambda}\xi^{p-1}a(t)+\cdots+{p\over\lambda^{p-1}}\xi
a(t)^{p-1}
\end{equation}
n'est pas
d'Eisenstein. En l'\'ecrivant
$$(a(t)\xi^{-1})^p=a(t)^pt^{-l}\Big(1+\cdots+{p\over\lambda^{p-1}}(a(t)\xi^{-1})^{p-1}\Big)
$$
nous allons pouvoir \'etablir le lemme analogue au lemme
\ref{aririri}, pourvu que
$a(T)\xi^{-1}\in B$,
pour $B$ la normalisation de $A[[t]]$ dans l'extension de Kummer
$L/K$. C'est le cas si nous supprimons le terme constant de $a(t)$ ; donc $a(t)t^{-l}$ est un polyn\^ome. Nous pouvons de
 nouveau utiliser le lemme \ref{lemBez} pour obtenir l'existence d'un
 \'el\'ement $\eta\in {\cal M}_B$ tel que $\xi=\eta^l$. Nous avons encore
\begin{lem}
\label{lemnorm1}
La normalisation de $A[[t]]$ dans l'extension de Kummer $L$ de $K$ est
l'alg\`ebre des s\'eries formelles $B=A[[\eta]]$.
\end{lem}
\begin{pr}
L'\'equation (\ref{equadefol}) montre que la norme de $\xi$ est $t^l$
; par suite la norme de $\eta$ est de la forme $\epsilon t$ avec
$\epsilon^l=1$, ce qui permet de supposer que la norme de $\eta$ est
$t$. Le polyn\^ome minimal de $\eta$ est de la forme
$F(t,\eta)=\eta^p+c_1\eta^{p-1}+\cdots+c_{p-1}\eta+t=0$
o\`u $c_i\in A[[T]]$. Comme dans le lemme \ref{lemnorm}, cette
\'equation implicite peut \^etre r\'esolue en $t$ ; autrement dit nous
pouvons exprimer $t$ comme s\'erie formelle en $\eta$. Donc $t\in
A[[\eta]]$. Alors l'alg\`ebre locale $A[[\eta]]$ qui est un
$A[[t]]$-module de type fini est de dimension de Krull
$\dim_{\Krull}A+1$, donc est l'alg\`ebre des s\'eries formelles en
$\eta$ sur $A$. Cette alg\`ebre est un sous-anneau int\'egralement
clos de $B$, donc $B=A[[\eta]]$.
\end{pr}
Analysons \`a pr\'esent cette d\'eformation au niveau
infinit\'esimal. Pour abr\'eger la discussion, nous nous limitons au cas
$l\not=1$, le cas $l=1$ se traitant de mani\`ere analogue. Fixons
$j\in[1,q-1]$ et sp\'ecialisons les variables
$x_i=0$ si $i\not=j$ et $x_j=\eps$ ($\eps^2=0$). Alors
$$a(t)=t^q+\eps t^{q-j}$$
La construction pr\'ec\'edente fournit une d\'eformation $\sigma$ de
l'automorphisme d'Artin-Schreier $\sigma_0$ \`a $k[\eps][[\eta]]$,
o\`u
$$\left\{\begin{array}{lll}\eta^l=\xi\\
t^l=\xi^p\Big(1-(a(t)\xi^{-1})^{p-1}\Big)\\
\sigma(\xi)=\xi+a(t)\\
\end{array}\right.$$
L'automorphisme initial $\sigma_0$ est d\'ecrit par les relations
$$\left\{\begin{array}{ll}\sigma_0(\xi)=\xi+t_0^q\\
t_0^l=\xi^p-t_0^{q(p-1)}\xi\\
\end{array}\right.$$
Regardons $t_0$ et $t$ comme fonctions implicites de la variable
$\eta$. \'Ecrivons
$$\sigma(\eta)=\sigma_0(\eta)+\eps\phi_j(\eta) \mbox{ avec } \phi_j\in
k[[\eta]]$$
D\'eterminons la direction tangentielle associ\'ee \`a cette
d\'eformation infinit\'esimale. Pour cela nous avons besoin du lemme
suivant :
\begin{lem} 
\label{lem425}
Nous avons
$v_{\eta}(\phi_j)=p(q-j)-(l-1)$.
\end{lem}
\begin{pr}
\'Ecrivons $t=t_0+\eps b$ pour $b\in k[[\eta]]$ et 
observons que $v_{\eta}(b)\geq 1$.
En effet
$$ t^l=t_0^l+\eps lt_0^{l-1}b=\xi^p-\xi t_0^{q(p-1)}+\eps(p-1)\xi
t_0^{(p-2)(q-j)}+\eps q(p-1)\xi bt_0^{q(p-1)-1}$$
D'o\`u $lb= (p-1)\xi t_0^{(p-2)(q-j)-(l-1)}+q(p-1)b\xi t_0^{q(p-1)-l}$
et $v_{\eta}(b)\geq 1$.
Par ailleurs
$$\sigma(\eta)^l=\sigma(\xi)=(\sigma_0(\eta)+\eps\phi_j(\eta))^l$$
Donc
$\sigma_0(\xi)+\eps l\sigma_0(\eta)^{l-1}\phi_j(\eta)=\xi+t^q+\eps
t^j$ et
$l\sigma_0^{l-1}(\eta)\phi_j(\eta)=qt_0^{q-1}b(\eta)+t_0^{q-j}$.
Ce qui donne la valuation de la s\'erie $\phi_j(\eta)$ :
$v_{\eta}(\phi_j)=p(q-j)-(l-1)$.
\end{pr} 
\begin{lem}
Les
classes
$\displaystyle\Big[\phi_j{d\over d\xi}\Big]\in H^1(G,\Theta)$
sont ind\'ependantes.
\end{lem}
\begin{pr}
Pour d\'emontrer ce lemme nous avons besoin d'une description de
$H^1(G,\theta)$ comme $k[[t_0]]$-module.
Pour cela, nous allons d'abord d\'efinir une $k[[t_0]]$-base de $k[[T]]$
relativement \`a laquelle $\delta$ est triangulaire sup\'erieur.
Pour tout $x,y\in k[[T]]$, par r\'ecurrence sur $n\geq 1$, il est
clair que 
$$\delta^n(xy)=\sum_{i,j\leq n,\; i+j\geq n}a_{ij}^n\delta^ix\delta^jy
\;\;\;\;\mbox{ avec } a_{i,n-i}=\binom{n}{i},\; 0\leq i\leq n$$
et que  
$$\forall\;\xi \;\mbox{ g\'en\'erateur de } E_2/E_1, \forall \;2\leq n\leq
p-1, \;\;\;\xi^n\in E_{n+1}\mbox{ et }\delta^n(\xi^n)=n!(\delta \xi)^n$$
Remarquons qu'il existe $\xi\in E_2-E_1$ tel que
si $0\leq j\leq i\leq p-1$ alors 
$$v_{T}(\delta^j(\xi^i))=pqj+l(i-j)
\mbox{ avec }l=v_{T}(\xi)\in[1,p[ \mbox{ et } m=pq-l$$
En effet,
nous pouvons choisir $\xi\in E_2$ dont la classe engendre le $k[[Y]]$-module $E_2/E_1$, nous pouvons m\^eme supposer que
$\xi\in Tk[[T]]$ ; alors pour $l=v_{T}(\xi)$ nous avons $1\leq l<p$ ;
car si $l\geq p$, alors $\xi Y^{-1}\in E_2$.
Comme $v_{T}(\delta \xi)=l+m$ et $\delta \xi\in E_1=k[[Y]]$, nous obtenons
$$v_{T}(\xi)=l\in[1,p[,\; v_{T}(\delta \xi)=l+m=pq$$
Ensuite nous raisonnons par r\'ecurrence sur $j$ en observant que si
$p$ ne divise pas $v_{T}(z)$ alors $v_{T}(\delta
z)=v_{T}(z)+m$. Posons :\\
$$\displaystyle\gamma(j)=\Big\lfloor{jpq+l(p-1-j)\over
  p}\Big\rfloor=jq+\Big\lfloor{l(p-1-j)\over p}\Big\rfloor$$
Alors 
$z_j=Y^{-\gamma(j)}\delta^j(\xi^{p-1}),$ $0\leq j\leq p-1$
 est un g\'en\'erateur de $E_{p-j}/E_{p-j-1}$.
En effet,
si $j=p-1$ alors $z_{p-1}=Y^{(p-1)l}(p-1)!(\delta \xi)^{p-1}$ est une unit\'e de $k[[Y]]$.
Si $j<p-1$, supposons que $z_j$ ne soit pas un g\'en\'erateur de $E_{p-j}/E_{p-j-1}$.
Il existe alors $\theta\in Y k[[Y]]$, $\xi\in E_{p-j}$ et $\eta\in E_{p-j-1}$ tels que
$z_j=\theta\xi+\eta$, donc $\delta^{p-j-1}(z_j)=\theta\delta^{p-j-1}(\xi)$.\\
Or d'une part $\delta^{p-j-1}(z_j)=Y^{-\gamma(j)}(p-1)!\delta(\xi)^{p-1}$ a pour valuation\\
$$v_{T}(\delta^{p-j-1}(z_j))=pq(p-1-j)-p\Big\lfloor{l(p-1-j)\over p}\Big\rfloor$$
D'autre part comme $v_{T}(\theta)\geq p$ et
$v_{T}(\delta^{p-j-1}(\xi))\geq v_{T}(\xi)+(p-j-1)m$, nous obtenons par 
comparaison
$$pq(p-1-j)-p\Big\lfloor{m(p-1-j)\over p}\Big\rfloor\geq p+(p-j-1)m$$ 
et vu
que $m=pq-l$,
$\displaystyle l(p-1-j)-p\Big\lfloor{l(p-1-j)\over p}\Big\rfloor\geq p$ ; ce qui est exclu.
Par cons\'equent
la famille $\{z_0,\ldots,z_{p-1}\}$ est une $k[[Y]]$-base de $k[[T]]$ avec
$$\delta(z_j)=Y^{\gamma(j+1)-\gamma(j)}z_{j+1}, 0\leq j\leq p-2
\mbox{ et } \delta(z_{p-1})=0
$$
Il ne reste plus qu'\`a voir que
la structure de $H^1(G,\Theta)$ comme $k[[Y]]$-module est donn\'ee par
$$H^1(G,\Theta)\cong \oplus_{j=1}^{p-1}{k[[Y]]\over
  (Y^{q+s_j})}$$
avec  $\displaystyle s_j\geq -1 \mbox{ et } s_{p-1}=\left\{\begin{array}{ll}0&\mbox{ si }
    l=1\\
-1&\mbox{ si } l\not=1\\
\end{array}\right.$.
 Nous observons d'abord que
$v_{T}(z_j)=r_j$ est le reste de
$l(p-1-j)$ modulo $p$ et que pour $0\leq j\leq p-1$, ces restes sont deux
\`a deux distincts.
Rappelons que $\Theta$ et la diff\'erente ${\cal D}$ sont
$G$-isomorphes et que 
${\cal D}=T^{\beta}k[[T]] \mbox{ avec } \beta=(m+1)(p-1)$.
Pour $(u_j)_{0\leq j\leq p-1}\in (k[[Y]])^p$,
on a 
$\sum_{j=0}^{p-1}u_jz_j\in{\cal D}$ si et seulement si 
pour tout $j$,
$$pv_{Y}(u_j)+r_j\geq (p-1)(pq-l+1)$$
c'est-\`a-dire si et seulement si
$$\left\{\begin{array}{ll}v_{Y}(u_j)\geq pq-(l-1)-q+1&\mbox{ si }r_j<l-1\\
v_{Y}(u_j)\geq pq-(l-1)-q&\mbox{ si }r_j\geq l-1\end{array}\right.$$
Donc 
$${\cal D}=\oplus_{j=0}^{p-1}k[[Y]]Y^{pq-(l-1)-q+s_j'}z_j \;\;\;\mbox
{ avec }\;\;\; s_j'=\left\{\begin{array}{ll}1 &\mbox{ si }r_j<l-1\\
0&\mbox{ si }r_j\geq l-1\\
\end{array}\right.$$
Posons $w_j=Y^{pq-(l-1)-q+s_j'}z_j$ pour $0\leq j\leq p-1$,
$s_j''=\gamma(j+1)-\gamma(j)-q$ et
$s_{j+1}=s_j'+s''_j-s_{j+1}'$ pour $0\leq j\leq p-2$. Alors
$\delta(w_j)=Y^{q+s_{j+1}}w_{j+1}$.
D'o\`u
$$H^1(G,\Theta)\cong \tors\Big({{\cal D}\over\delta({\cal
    D})}\Big)\cong\oplus_{j=1}^{p-1}{k[[Y]]\over (Y^{q+s_j})}$$
Nous pouvons maintenant d\'emontrer le lemme annonc\'e.\\
En effet par le lemme \ref{lem425},
$$v_{\eta}(\phi_j{dt_0\over d\xi})=p(q-j)-(l-1)+\beta=p(pq-(l-1)-q+q-j)$$
Par cons\'equent la composante de $\displaystyle \phi_j{dt_0\over d\xi}$ sur
$w_{p-1}=t_0^{pq-(l-1)-q+s'_{p-1}}z_{p-1}$ est non nulle. Le
coefficient relatif \`a $w_{p-1}$ de $\displaystyle \phi_j{dt_0\over d\xi}$
a pour valuation en $t_0$, $q-j-s'_{p-1}$. Donc ces composantes dans
$H^1(G,\Theta)^G$ sont ind\'ependantes.
\end{pr}
\begin{thm}
\label{result}
Soit $m=pq-l$, $q\geq 1$ et $l\in[1,p-1]$. L'anneau versel pour un automorphisme
$\sigma_0$ de conducteur $m$ a un quotient formellement lisse :
$$\left\{\begin{array}{ll}R_{\sigma_0}\twoheadrightarrow
   W(k)(\zeta)[[x_1,\ldots,x_q]]&\mbox{ si } l=1\\
R_{\sigma_0}\twoheadrightarrow
    W(k)(\zeta)[[x_1,\ldots,x_{q-1}]]&\mbox{ si } l\not=1\\
\end{array}\right.$$
Si $p=2$, pour tout conducteur $m$ (impair), l'anneau de d\'eformations
versel d'un automorphisme d'ordre 2 est l'alg\`ebre des s\'eries formelles
$$W(k)[[X_1,\ldots,X_{{m+1\over 2}}]]$$
Il n'y a donc pas d'obstruction au probl\`eme de rel\`evement.
\end{thm}
\begin{rem}
Le th\'eor\`eme \ref{result} a pour cons\'equence g\'eom\'etrique en caract\'eristique $2$, le
fait qu'un point de branchement de conducteur $m$ se d\'eploie dans
une d\'eformation g\'en\'erique en
$\displaystyle d={m+1\over 2}$ points de conducteur 1. Notons aussi que
le th\'eor\`eme \ref{result} joint au th\'eor\`eme \ref{thmloglo} 
redonne le r\'esultat de Laudal et L\o nsted ([LaL\o] Theorem 2), c'est-\`a-dire
si $C$ est hyperelliptique et $\sigma$ est l'involution
hyperelliptique alors l'anneau (uni)versel de d\'eformations global est
une $W(k)$-alg\`ebre de s\'eries formelles.\\
Dans le cas g\'en\'eral ($G$ cyclique d'ordre $p$), le th\'eor\`eme
\ref{result} entra\^\i ne une minoration de la dimension de Krull de
l'anneau versel local. Dans le paragraphee suivant, nous montrons que
cette minoration est une \'egalit\'e.
\end{rem}
\section{Un espace de modules pour les rev\^etements galoisiens}
\label{subsec01}
Fixons $n\geq 3$ premier \`a $p$.
Rappelons la d\'efinition d'une structure de niveau $n$
(\cite{DeMu}, \cite{KaMa}). Soit $C$ une courbe lisse projective de
genre $g\geq 1$  et $J(C)$ la jacobienne de $C$.
Notons $C[n]=J(C)[n]$ le sous-groupe des points de $n$-division de
$J(C)$. Rappelons que 
pour $S$ un $\Z[1/n]$-sch\'ema et $C$ une $S$-courbe
projective et lisse, le $S$-sch\'ema en groupe
$$C[n]=J(C)[n]\rightarrow S$$
est fini \'etale sur
$S$ (\cite{KaMa} Proposition 1.6.4).

Une {\it structure de niveau $n$}\index{structure de niveau $n$} pour la courbe $C/S$,
est un isomorphisme $\varphi$ de $S$-sch\'emas en groupes
$$
C[n]\stackrel{\varphi\sim}
  \longrightarrow(\Z/n\Z)_S^{2g}=(\Z/n\Z)^{2g}\times S$$
Pour $S$ un sch\'ema tel que $n$ soit inversible dans ${\cal O}_S$,
notons $(C/S,\varphi)$ le couple compos\'e d'une $S$-courbe
$\pi:C\rightarrow S$ lisse et projective, et d'une
structure de niveau $\varphi$ (le niveau est fix\'e \'egal \`a $n$). Un
isomorphisme $\tau:(C/S,\varphi)\stackrel\sim\longrightarrow (C'/S,\varphi')$
est un $S$-isomorphisme $\tau:C\stackrel\sim\longrightarrow C'$ tel que le triangle
ci-dessous soit commutatif :
\begin{diagram}[labelstyle=\textstyle]
f:  C'[n]&&\rTo^{\tau[n]}_{\sim}&&C[n]\\
&\rdTo^{\sim}_{\varphi'}&&\ldTo_{\varphi}^{{\sim}}\\
&&(\Z/n\Z)_S^{2g}
\end{diagram}
o\`u $\tau[n]$ est l'isomorphisme induit par l'image r\'eciproque
${\cal L}\rightarrow\tau^*({\cal L})$ de $\Pic(C')$ sur $\Pic(C)$.
Rappelons la propri\'et\'e de rigidit\'e (\cite{KaMa}):
\begin{propf}
\label{thrigi}{\it
Si $S$ est connexe, il y a \'equivalence entre les propri\'et\'es
suivantes :\\
i. $\tau[n]$ est l'identit\'e en un point $s\in S$\\
ii. $\tau[n]$ est l'identit\'e\\
iii. $\tau$ est l'identit\'e.}
\end{propf}
Soit $\Sch_{\Z[1/n]}$ la cat\'egorie des sch\'emas $S$ tels que $n$ est
  inversible sur ${\cal O}_S$ ; d\'efinissons un foncteur
  contravariant
${\cal M}_{g,n}$ ($g\geq 1)$ par
$${\cal M}_{g,n}:\Sch_{\Z[1/n]}\rightarrow \Ens,\;\;\;\;\; S\mapsto
\Big\{\mbox{classes d'\'equivalence }  [C/S,\varphi]\Big\}$$
Rappelons le r\'esultat fondamental (\cite{DeMu}, \cite{KaMa} pour
$g=1$) :
\begin{thmf}{\it 
Le foncteur ${\cal M}_{g,n}$ est
repr\'esentable par un $\Z[1/n]$-sch\'ema quasi-projectif not\'e
$\M_{g,n}$ appel\'e le sch\'ema modulaire des courbes de genre $g$ de
niveau $n$.}
\end{thmf}
La repr\'esentabilit\'e de ${\cal M}_{g,n}$  signifie en particulier
qu'il existe une $\M_{g,n}$-courbe universelle avec
structure de niveau $n$
$$ {\cal C}_{g,n}\rightarrow \M_{g,n}$$
Si $(C/S,\varphi)$ est une $S$-courbe projective lisse munie d'une
structure de niveau $n$, il existe un unique morphisme $f:S\rightarrow
\M_{g,n}$ tel que $C$ soit $S$-isomorphe \`a ${\cal
  C}_{g,n}\times_{\M_{g,n}} S$ et $\varphi$ s'identifie \`a l'image
  r\'eciproque de la structure de niveau universelle sur ${\cal C}_{g,n}$.
Si $S$ est connexe,
le groupe $\Gl_{2g}(\Z/n\Z)_S$ est le groupe des automorphismes du
$S$-sch\'ema constant $(\Z/n\Z)_S^{2g}$. Ainsi $\Gl_{2g}(\Z/n\Z)_S$ op\`ere sur l'ensemble des structures
de niveau $n$ port\'ees par la courbe $C/S$ par :
$$\alpha\in\Gl_{2g}(\Z/n\Z)_S,\;\;\;\;
\alpha[C/S,\varphi]=[C/S,\alpha\circ\varphi]$$
Dire que $\alpha[C/S,\varphi]=[C/S,\varphi]$ signifie qu'il existe un
automorphisme $\tau$ (unique par la propri\'et\'e de rigidit\'e) de
$C/S$ tel que le diagramme suivant soit commutatif 
\begin{diagram}[labelstyle=\textstyle]
C[n]&\lTo_{\tau[n]}^{ }&C[n]\\
\dTo_{\varphi}^{\sim}&&\dTo^{\sim}_{\varphi}\\
(\Z/n\Z)_S^{2g}&\lTo^{\alpha}&(\Z/n\Z)_S^{2g}\\
\end{diagram}
Soit $G$ un sous-groupe
(fini) du groupe $\Aut(C/S)$ et soit $\varphi$ une structure de niveau
$n$ sur $C/S$.
\begin{lemf}{\it
Supposons $S$ connexe. Pour tout point $s\in S$, l'image $G_s$
de $G$  dans $\Gl_{2g}(\Z/n\Z)$
$$G\hookrightarrow\Aut(C_s)\rightarrow\Aut(C_s[n])\stackrel{\varphi_s \sim}\longrightarrow
\Gl_{2g}(\Z/n\Z)$$
est  ind\'ependante de $s\in S$; autrement dit,
$\varphi(G)$ est un sous-groupe constant de $\Gl_{2g}(\Z/n\Z)_S$.}
\end{lemf}
Fixons \`a pr\'esent un sous-groupe $G\subset \Gl_{2g}(\Z/n\Z)$. Notons
$(C/S,\varphi,\theta)$ le triplet constitu\'e
d'une $S$-courbe alg\'ebrique lisse compl\`ete de genre $g\geq 1$ 
munie d'une struture $\varphi$ de niveau $n$ et d'une action
fid\`ele $\theta$ de $G$ sur $C/S$ qui est $\varphi$-adapt\'ee,
c'est-\`a-dire telle que l'isomorphisme $\varphi$ soit
$G$-\'equivariant ; ou encore que pour tout $\sigma\in G$
$$\varphi\circ\theta(\sigma)^{-1}[n]=\sigma\circ\varphi$$
Deux triplets $(C/S,\varphi,\theta)$ et $(C'/S,\varphi',\theta')$ sont
isomorphes s'il existe un $S$-isomorphisme
$\tau:C\stackrel\sim\longrightarrow C'$ qui est $G$-\'equivariant et
qui rend le diagramme suivant commutatif
\begin{diagram}[labelstyle=\textstyle]
C[n]&&\lTo^{\tau[n]}
  &&C'[n]\\
&\rdTo{\sim }{\varphi}&&\ldTo{\sim}{\varphi'}\\ \;\;\;\;
&&(\Z/n\Z)^{2g}_S\\
\end{diagram}
Notons que si $\tau[n]$ est $G$-\'equivariant, alors $\tau$ est
$G$-\'equivariant d'apr\`es le th\'eor\`eme de rigidit\'e \ref{thrigi}.
Si $n$ est premier \`a $p$, nous pouvons d\'efinir le
foncteur contravariant ${\cal M}_{g,n,G}$ par
$${\cal M}_{g,n,G}:S\mapsto \{
\mbox{classes d'\'equivalence } 
[C/S,\varphi,\theta]\}$$ 
\begin{thmf}{\it
Le foncteur ${\cal M}_{g,n,G}$ 
est repr\'esentable par le
sous-sch\'ema ferm\'e des points fixes de $G$ sur $\M_{g,n}$ :
$$\M_{g,n,G}=\M_{g,n}^G$$}
\end{thmf}
\begin{pr}
Observons tout d'abord que de l'action naturelle de $\Gl_{2g}(\Z/n\Z)$
sur $\M_{g,n}$ d\'ecoule une action de $G$. Le
sous-sch\'ema des points fixes de $G$ a pour points
$$\Hom(S,\M_{g,n}^G)\cong\Hom(S,\M_{g,n})^G$$
autrement dit le morphisme $f:S\rightarrow \M_{g,n}$ se factorise par
$\M_{g,n}$ si et seulement si pour tout $g\in G$, $g\cdot f=f$.
Il s'agit alors d'identifier les points de ${\cal M}_{g,n,G}(S)$ avec
ceux de ${\cal M}_{g,n}(S)^G.$\\
Soit $[C/S,\varphi,\theta]\in {\cal M}_{g,n,G}(S)$, par le foncteur d'oubli
de l'action de $G$, on obtient un couple $(C/S,\varphi)$ dont la classe
\`a isomorphisme pr\`es est bien d\'efinie ; on a donc un point
$[C/S,\varphi]\in {\cal M}_{g,n}(S)$. Ce point est fix\'e par $G$, car,
par d\'efinition, pour tout $g\in G$, $\theta(g)\in \Aut(C/S)$ 
induit exactement $g$ via $\varphi$.\\
Ce foncteur d'oubli est bien injectif. En effet si
$[C/S,\varphi,\theta]$ et $[C'/S,\varphi',\theta']$ ont m\^eme image
alors $[C/S,\varphi]=[C'/S,\varphi']$. Donc il existe un
$S$-isomorphisme 
$$\tau:C\stackrel\sim\rightarrow C'$$
qui rend le diagramme suivant commutatif
\begin{diagram}[labelstyle=\textstyle]
C[n]&&\lTo^{\tau[n]}_{\sim}
  &&C'[n]\\
&\rdTo{\sim}{\varphi}&&\ldTo{\sim}{\varphi'}\\ \;\;\;\;
&&(\Z/n\Z)^{2g}_S\\
\end{diagram} 
Les actions de $G$ via $\varphi$ et $\varphi'$ sont
$G$-\'equivariantes, donc $\tau[n]$ est $G$-\'equivariant ; par
rigidit\'e, $\tau$ est aussi $G$-\'equivariant. Donc $\tau$ est un
isomorphisme entre les triplets $(C/S,\varphi,\theta)$ et
$(C'/S,\varphi',\theta')$.\\
Pour montrer la surjectivit\'e, prenons $[C/S,\varphi]\in{\cal
  M}_{g,n}(S)^G$. Par d\'efinition pour tout $g\in G$ il existe un
unique automorphisme $\tau(g)$ :
$$\tau(g):C\stackrel\sim\longrightarrow C$$
qui est $\varphi$-compatible, c'est-\`a-dire tel que le diagramme
suivant commute
\begin{diagram}[labelstyle=\textstyle]
C[n]&&\lTo^{\tau(g)[n]}_{\sim}
  &&C[n]\\
&\rdTo{\sim}{\varphi}&&\ldTo{\sim}{\varphi}\\ \;\;\;\;
&&(\Z/n\Z)^{2g}_S\\
\end{diagram} 
Ensuite on observe que le morphisme $\theta : g\mapsto\tau(g)$ d\'efinit une
action de $G$ sur $C/S$. D'o\`u la surjectivit\'e. 
\end{pr}
Il existe une courbe universelle avec structure de niveau $n$ et
action de $G$ d\'efinie sur $\M_{g,n,G}$ :
$$\pi_{g,n,G} : {\cal C}_{g,n,G}\rightarrow \M_{g,n,G}$$
Il se peut que $\M_{g,n,G}$ ne soit pas connexe. Alors  si pour tout
$x\in \M_{g,n,G}$, nous notons ${\cal C}_{x}=\pi_{g,n,G}^{-1}(x)$, le
genre de ${\cal C}_{x}/G$ n'est pas forc\'ement constant. Nous fixons alors un
genre $g'$ et nous nous limitons \`a l'\'etude des seules composantes
connexes de $\M_{g,n,G}$ le long desquelles le genre des courbes
${\cal C}_{x}/G$ est constant et \'egal \`a $g'$ (nous supposons ce probl\`eme
non vide!). 
\subsection{L'espace modulaire de Harbater}
Fixons une courbe affine lisse $U$. Soit $\Sigma$ sa compl\'etion et
$G$ un $p$-groupe fini. Harbater (\cite{Har}) a \'etudi\'e la
vari\'et\'e des $G$-rev\^etements \'etales de $U$, c'est-\`a-dire de
$\Sigma$ avec points de branchements contenus dans $\Sigma-U$. Ici, et
contrairement \`a la situation consid\'er\'ee dans les paragraphes
pr\'ec\'edents, deux rev\^etements galoisiens \'etales de $U$ de
groupe $G$, $\pi:C\rightarrow U$ et $\pi':C'\rightarrow U$ sont
\'equivalents s'il existe un isomorphisme
$\tau:C\stackrel\sim\longrightarrow C'$ tel que $\pi'\tau=\pi.$ Dans
\cite{Har}, Harbater prouve que les classes de $G$-rev\^etements
\'etales de $U$ sont les points d'un espace de modules fin, qui est
une limite inductive d'espaces affines.\\

Nous allons d\'etailler le cas particulier $G=\Z/p\Z$ et
$U=\P^1-\{b_1,\ldots,b_r\}$. Dans ce cas un rev\^etement galoisien
\'etale de $U$, c'est-\`a-dire une courbe $C$, est d\'ecrit par le
corps $k(C)=k(t)(\xi)$, avec $\xi$ solution de l'\'equation
d'Artin-Schreier
$$\xi^p-\xi=a(t)\in\Gamma(U,{\cal O}_U)$$
La classe du rev\^etement est d\'etermin\'ee par la classe de la
fraction rationnelle $a(t)$ modulo l'image de $\Gamma(U,{\cal O}_U)$
par l'application $\gp:x\mapsto x^p-x$. D'o\`u la suite exacte qui
d\'ecrit l'ensemble des points du sch\'ema de Harbater $\H(U)$ (\cite{Har})
\begin{equation}
\label{suiHar}
0\longrightarrow {\Gamma(U,{\cal O}_U)\over k}\stackrel{\gp}\longrightarrow {\Gamma(U,{\cal O}_U)\over k}\longrightarrow
\left\{\begin{array}{ll}\mbox{classes de}\cr
\mbox{rev\^etements}\end{array}\right\}\rightarrow 0
\end{equation}
Notons $[a(t)]\in H(U)$ l'image de $a(t)\in\Gamma(U,{\cal O}_U)$.\\
Soit $b\in \Sigma-U$, un point de branchement ; soit $t'$ le
param\`etre local $t-t(b)$ (resp. $t^{-1}$ si $b=\infty$) en $b$. De
mani\`ere analogue \`a (\ref{suiHar}), les extensions
$\Z/p\Z$-cycliques de $k((t))$ sont classifi\'ees par un
g\'en\'erateur d'Artin-Schreier, c'est-\`a-dire, par la classe de $a\in
k((t'))$ modulo $k[[t']]+\gp(k((t')))$.\\
Nous avons l'analogue local de (\ref{suiHar}),  
\begin{equation}
\label{suiHar2}
0\longrightarrow {k((t'))\over
  k[[t']]}\stackrel\gp\longrightarrow{k((t'))\over
  k[[t']]}\longrightarrow\left\{\begin{array}{ll}\mbox{classes de rev\^etements}\cr
\mbox{de } \Spec k((t'))\end{array}\right\}\rightarrow 0
\end{equation}
\begin{dea}
Appellons {\it partie
  polaire} en $b$,
un polyn\^ome de Laurent
$$\gt=\sum_{j=1}^m{\alpha_j\over {t'}^j},\;\;\mbox{ o\`u }\;\;
\alpha_j=0 \mbox{ si } p|j$$
et $\alpha_m\not=0$. L'entier $m$, premier \`a $p$, est la {\it longueur}\index{longueur} de $\gt$. 
\end{dea}
\begin{lem}{\it
Dans chaque classe de $$\displaystyle {k((t'))\over
  k[[t']]+\gp(k((t')))}$$ il y a une unique partie polaire.}
\end{lem}
\begin{pr}
Consid\'erons la classe de $a(t')$ pour 
$$a(t')=\sum_{j=1}^m{\alpha_j\over{t'}^j}\in k((t'))$$
Montrons d'abord qu'il existe une partie polaire dans la classe de
$a$. Supposons que $a$ n'est pas une partie polaire.
Soit $j(a)$ le plus grand indice $j\in[1,m]$ tel que $p$ divise $j$ et
$\alpha_j\not=0$. Notons $j(a)=pl$ et $\alpha_{j(a)}=\beta^p$ pour
$\beta\in k^*$. Alors 
$$b(t')=a(t')-\gp\Big({\beta\over{t'}^l}\Big)$$
appartient \`a la classe de $a$ et si $b$ n'est pas une partie polaire,
$j(b)<j(a)$. Il est alors imm\'ediat que la classe de $a$ contient une
partie polaire et que celle-ci est unique.
\end{pr}
Soit maintenant, $\pi:C\rightarrow\P^1$ un $G$-rev\^etement \'etale
sur $U=\P^1-\{b_1,\ldots,b_r\}$, d\'etermin\'e par la classe de
$a(t)\in\Gamma(U,{\cal O}_U)$. Pour tout point $b_i\in \P^1-U$,
$1\leq i\leq r$,
soit $\gt_i$ la partie polaire de $a(t)$ dans $k((t'_i))$ et soit
$m_i$ sa longueur. Remarquons que $m_i$ est le conducteur du
rev\^etement au point $b_i$. Nous obtenons ainsi une application
d\'efinie
sur l'espace modulaire de Harbater $\H(U)$ de $U$ :
\begin{equation}
\label{appliHar}
[a(t)]\mapsto
\{\gt_1,\ldots,\gt_r\}
\end{equation}
Nous avons (\cite{Har} Proposition 2.7) :
\begin{prop}{\it
L'application (\ref{appliHar}) est un isomorphisme de l'espace
modulaire de Harbater $\H(U)$ sur les $r$-uplets de rev\^etements
locaux aux points $b_i\in \P^1-U$, $1\leq i\leq r$.}
\end{prop}
La relation entre $a(t)$ et les parties polaires
$\{\gt_1,\ldots,\gt_r\}$ est la d\'ecomposition en
\'el\'ements simples
$$a(t)=\sum_{i=1}^r\gt_i$$
Fixons \`a pr\'esent en chaque point $b_i$, $1\leq i\leq r$, le conducteur $m_i$
(c'est-\`a-dire la longueur de la partie polaire $\gt_i$). L'espace
modulaire correspondant $\H(U,\{m_i\}_{1\leq i\leq r})$ v\'erifie
$$ \H(U,\{m_i\}_{1\leq i\leq r})\cong (\A_*^1)^r\times(\A^1)^{r'}$$
avec $r'=\sum_{i=1}^r (m_i-\lfloor m_i/p\rfloor -1)$. En
particulier sa dimension est
$$\dim_k\H(U,\{m_i\}_{1\leq i\leq
  r})=\sum_{i=1}^r\Big(m_i-\Big\lfloor{m_i\over p}\Big\rfloor\Big)$$
Citons aussi le corollaire suivant du r\'esultat de Harbater dont nous
ferons usage pour $G=\Z/p\Z$ (\cite{Har} Corollary 2.4) :
\begin{cor}{\it
Soit $G$ un $p$-groupe fini et $U=\A^1$. Tout $G$-rev\^etement de $\Spec k((t'))$ ($t'=t^{-1}$) se prolonge de mani\`ere
unique en un $G$-rev\^etement de $\P^1$, \'etale sur
$\A^1=\P^1-\{\infty\}$.}
\end{cor} 
\subsection{D\'etermination de la dimension de Krull}
Soit $R_{\sigma}$ l'anneau
de d\'eformations versel local associ\'e \`a un automorphisme $\sigma$ d'ordre
$p$ de conducteur $m$ de $k[[T]]$. Le groupe cyclique d'ordre $p$, $G=<\sigma>$
s'identifie au groupe d'inertie \`a l'infini d'un $G$-rev\^etement
$\pi:C\rightarrow \P^1$ \'etale sur $\A^1$. Fixons une structure de
niveau $n$ sur la courbe $C$. Posons $N=m+1$. Le genre $g$ de $C$ est
donn\'e par la relation de Riemann-Hurwitz :
$$g={(N-2)(p-1)\over 2}$$
Nous pouvons exclure le cas $m=1$ car alors l'anneau $R_{\sigma}$ est
d\'ecrit explicitement (Th\'eor\`eme \ref{resultcond}) et nous pouvons
\'egalement supposer $p\not=2$ (Th\'eor\`eme
\ref{result}). Dans la suite nous nous pla\c{c}ons donc dans le cas
$p>2$ et $m>1$. Ainsi le genre du rev\^etement $C$ de $\P^1$, \'etale
sur $\A^1$, est $g\geq 2$ sauf si $(p,m)=(3,2)$. Dans cette situation
exceptionnelle, nous pouvons remplacer $C$ par un rev\^etement de
$\P^1$ ramifi\'e en deux points $0$ et $\infty$, avec le conducteur
$m=3$ en chaque point. Donc dor\'enavant les rev\^etements
consid\'er\'es sont de genre $g\geq 2$.
Un tel rev\^etement $\pi:C\rightarrow\P^1$, d\'efini sur $k$, se r\'ealise comme un point $x$
d'un espace modulaire $\M=\M_{g,n,G}$. Cet espace est d\'efini sur
$\Z[1/n]$. Par changement de l'anneau de base, nous pouvons supposer 
qu'il est d\'efini sur $W(k)$. Le point $x$ appartient alors \`a la
fibre sp\'eciale de $\M$.  
\begin{lem}{\it
L'anneau $R_{\glo}$ de d\'eformations
(uni)versel de $(C,G)$ est reli\'e \`a $\M$ par
$$R_{\glo}\cong\widehat{\cal O}_{\M,x}$$}
\end{lem}
\begin{pr}
L'anneau local complet $\widehat{{\cal O}}_{\M,x}$ pro-repr\'esente le foncteur
$$F:{\cal C}\rightarrow \Ens,\;\;\; A\mapsto\left\{\begin{array}{ll}
\mbox{rel\`evement } \alpha\in \Hom(\Spec A,\M) \cr
\mbox{ de } x\in\Hom(\Spec k,\M)\cr
\end{array}\right\}$$
Par oubli de la structure de niveau, nous d\'efinissons un morphisme lisse
$$\phi:F\rightarrow D_{\glo}$$
Et 
$$d\phi:F(k[\eps])\stackrel\sim\longrightarrow D_{\glo}(k[\eps])$$
est bijectif. En effet si $(X,G)$ est une d\'eformation de $(C,G)$ \`a
$k[\eps]$, alors la structure de niveau fix\'ee sur $C$ se rel\`eve de
mani\`ere unique \`a $X$. Par d\'efinition, nous obtenons
$R_{\glo}\cong\widehat{\cal O}_{\M,x}$.
\end{pr}
La dimension de Krull de $\widehat{\cal O}_{\M,x}$, c'est-\`a-dire la
dimension de $\M$, est connue lorsque
$p$ ne divise pas l'ordre de $G$ (ramification mod\'er\'ee). 
Nous allons prouver que cette dimension est ind\'ependante de
l'hypoth\`ese sur l'ordre de $G$, c'est-\`a-dire m\^eme si $G$ est
cyclique d'ordre $p$.
\begin{prop}
\label{propf33}
{\it
Nous avons
$\dim_{\Krull}\widehat{{\cal O}}_{\M,x}=N-2$.}
\end{prop}
Comme dans le cas de la ramification mod\'er\'ee, nous souhaitons comparer
$\M$ \`a un sch\'ema de ``configuration'', celui form\'e par les
points de branchement, suppos\'es mobiles. Le probl\`eme est que ce
nombre n'est pas constant dans une d\'eformation, contrairement au cas
mod\'er\'e. Nous contournons cette difficult\'e de la mani\`ere
suivante (voir \cite{Be} \S 2 Proposition) :
\begin{lem}
\label{lemf33}
{\it
Soit $\pi:X\rightarrow S$ un $S$-point de $\M$, $X$ et $S$
$k$-sch\'emas  de type fini ($k$ alg\'ebriquement clos). Le sch\'ema des points
fixes $\Fix=X^G$ de $G$ (cyclique) sur $X$ est un diviseur de Cartier relatif de
degr\'e
$N$.}
\end{lem}
\begin{pr}
Soit $s\in S$ et $x\in X_s=\pi^{-1}(s)$ un point de $\Fix=X^G$ ; $x$ est un point
ferm\'e de la fibre $X_s$ car l'automorphisme $\sigma$ n'est pas l'identit\'e
sur $X_s$. L'anneau local ${\cal O}_{X_s,x}$ est de valuation
discr\`ete. Choisissons une uniformisante $t$ de ${\cal O}_{X_s,x}$ que nous relevons en
$T\in{\cal M}_x\subset{\cal O}_{X,x}$. Nous avons alors
$${\cal M}_x={\cal M}_s{\cal O}_{X,x}+T{\cal O}_{X,x}$$
Nous allons prouver que l'id\'eal ${\cal I}$ du sous-sch\'ema $\Fix$ a
pour \'equation locale en $x$, $\sigma^*(T)-T$. Un \'el\'ement $y\in
{\cal M}_x$ peut se mettre sous la forme 
$$y\equiv P_l(T) \mbox{ mod } {\cal M}_x^{l+1}$$
o\`u $l\geq 1$ et $P_l(T)$ est un polyn\^ome \`a coefficients dans
l'image par $\pi^*$ de ${\cal O}_{S,s}$  dans ${\cal O}_{X,x}$ de degr\'e inf\'erieur \`a $l$ ; autrement dit
$$y\equiv\sum_{i\leq l}r_iT^i \mbox{ mod } {\cal M}_x^{l+1}$$
Alors modulo ${\cal M}_x^{l+1}$,
$$ \sigma^*(y)-y=\sum_{i\leq l}r_i(\sigma^*(T^i)-T^i)=\sum_{i\leq
  l}r_i(\sigma^*(T)^i -T^i)$$
Or $\sigma^*(T)^i-T^i=(\sigma^*(T)-T)\lambda_i$ avec
$\lambda_i\in{\cal O}_{X,x}$, ainsi
$$\sigma^*(y)-y\in{\cal O}_{X,x}(\sigma^*(T)-T)+{\cal M}_x^{l+1}$$ et
${\cal I}_x\subset{\cal O}_{X,x}(\sigma^*(T)-T)+{\cal M}_x^{l+1}$
pour tout $l\geq 1$, ce qui entra\^\i ne l'\'egalit\'e
$${\cal I}_x={\cal O}_{X,x}(\sigma^*(T)-T)$$
Il est clair que $\Fix$ est un diviseur de Cartier relatif. La fibre
de $\Fix$ en un point g\'eom\'etrique $s$ de $S$ est le sous-sch\'ema
des points fixes de $\sigma$ sur la fibre $X_s$, donc
$$\Fix_s=\sum_{i=1}^{r'}(a_i+1)w_i$$
o\`u $w_1,\ldots,w_{r'}$ sont les points fixes de conducteurs respectifs
$a_1,\ldots,a_{r'}$ ($a_i=0$ si $p\not=0$ dans $k(s)$). Le degr\'e est
$N=\sum_{i=1}^{r'}(a_i+1)$.
\end{pr}
Le lemme \ref{lemf33} nous permet de d\'emontrer la proposition
\ref{propf33}, de la fa\c{c}on suivante :
D'apr\`es le lemme \ref{lemf33} pour la courbe
universelle
$${\cal C}_{g,n,G}\rightarrow\M$$
nous obtenons un diviseur de Cartier relatif not\'e ${\cal F}={\cal
  F}_{g,n,G}\subset {\cal C}_{g,n,G}$.\\
Soit $\Sigma={\cal C}_{g,n,G}/G$ et $\tau:{\cal
  C}_{g,n,G}\rightarrow\Sigma$, $h:\Sigma\rightarrow\M$ les morphismes
correspondants. Remarquons que $\tau$ est fini et plat de rang
$p=|G|$. Il est alors possible de former l'image directe $\tau_*({\cal
  F})={\cal D}$, qui est un diviseur de Cartier relatif de $\Sigma$
au-dessus de $\M$ de degr\'e $N$ (\cite{Mu} Lecture 10). Par
hypoth\`ese, les fibres de $h$ sont de genre z\'ero, c'est-\`a-dire que
les fibres g\'eom\'etriques sont des droites projectives ($h$ est un
fibr\'e tordu en droites projectives).\\

Notons $\M_0$ la fibre sp\'eciale de
 $\M\rightarrow
\Spec(W(k))$. Si $Z\subset \M_0$ est une composante irr\'eductible
avec $x\in Z$ et $Z$ horizontale ($Z\not\subset M_0$), alors $\dim
Z=1+\dim Z_{\eta}$ o\`u $Z_{\eta}$ est la fibre g\'en\'erique. 
Du fait que les rev\^etements cycliques de degr\'e $p$ se rel\`event
en caract\'eristique z\'ero, une telle composante horizontale existe. 
Dans le
cas de la ramification mod\'er\'ee, la dimension de $\M$ est connue,
nous avons $\dim Z_{\eta}=N-3$, donc $\dim Z=N-2$.\\

 Nous nous limitons
\`a pr\'esent aux composantes verticales. Pour estimer $\dim {\cal
  O}_{\M_0,x}$, nous pouvons consid\'erer un voisinage \'etale de $x$
dans $\M_0$. En particulier comme $h:\Sigma\rightarrow \M$ est lisse,
$h$ poss\`ede, localement pour la topologie \'etale, une section. Nous
pouvons donc supposer que la courbe induite $\Sigma_V\rightarrow V$
est un produit $V\times\P^1$. Fixons une trivialisation
$\phi:\Sigma_V\stackrel\sim\longrightarrow V\times \P^1$.
Cela signifie que pour tout point g\'eom\'etrique $v\in V$, nous avons
une identification privil\'egi\'ee
$$\phi_v:\Sigma_v={\cal C}_v/G\stackrel{\sim}\longrightarrow \P^1$$
La fibre ${\cal D}_v$ de ${\cal D}$ est un diviseur de $\Sigma_v$ :
$${\cal D}_v=\sum_{i=1}^{r'}(a_i+1)b_i$$ 
o\`u $b_1,\ldots,b_{r'}$ sont les points de branchement de ${\cal
  C}_v\rightarrow\Sigma_v\cong\P^1$ et $a_1,\ldots,a_{r'}$ les
conducteurs respectifs. Nous d\'efinissons ainsi un morphisme
$$\delta:V\rightarrow \Div_N(\P^1),\;\;\;\; v\mapsto\phi_v({\cal
  D}_v)$$
Soit $\PGL(1)$ le groupe des automorphismes de la droite projective
$\P^1_k$. Formons le morphisme
$$\Delta:V\times\PGL(1)\rightarrow\Div_N(\P^1),\;\;\;\;
\Delta(v,h)=h(\delta(v))$$
(Nous utilisons la m\^eme notation $h$ pour un \'el\'ement de
$\PGL(1)$ et son action sur $\Div_N(\P^1)$). Nous pouvons supposer
que 
$$\Delta(x,1)=D_0=(m+1)\infty$$
Soit $D=\sum_{i=1}^{r'}(a_i+1)b_i$ un point g\'eom\'etrique de
$\Div_N(\P^1)$. En particulier $\sum_{i=1}^{r'}(a_i+1)=N+1$. Examinons la
fibre $\Delta^{-1}(D)$. Un point de cette fibre est un couple
$(z,h)\in V\times \PGL(1)$ tel que $h(\phi_v(D_v))=D$. Si nous
rempla\c{c}ons les courbes $\Sigma_V$ et ${\cal C}_V$ par leurs images
r\'eciproques par la projection $V\times\PGL(1)\rightarrow V$, nous
obtenons par restriction \`a $Z=\Delta^{-1}(D)$ le diagramme
 \begin{diagram}[labelstyle=\textstyle]
{\cal C}_Z&&\rTo  &&Z\\
&\rdTo&&\ldTo\\ \;\;\;\;
&&\Sigma_Z={\cal C}_Z/G\\
\end{diagram} 
En corrigeant l'identification
$\phi_v:\Sigma_{(v,h)}\stackrel\sim\rightarrow\P^1$ par
$h\circ\phi_v$, nous obtenons un isomorphisme $\Sigma_Z\cong Z\times
\P^1$, tel que pour tout $z\in Z$, le diviseur de branchement de
${\cal C}_z\rightarrow\Sigma_z=\P^1$ est constant \'egal \`a $D$. Par
restriction \`a $U=\P^1-\Supp(D)$, nous obtenons une famille de
$G$-rev\^etements \'etales de $U$ param\'etr\'ee par $Z$ ; d'o\`u un
morphisme
$Z\rightarrow \H(U,\{a_i\}_{1\leq i\leq r'})$ dans l'espace modulaire de Harbater des
$G$-rev\^etements \'etales de $U$ avec conducteurs fix\'es en les
points de $\P^1-U$. Ce morphisme est \`a fibres finies, donc
$$\dim Z\leq\dim_k\H(U,\{a_i\}_{1\leq i\leq r'})$$
soit encore
$$\dim Z\leq \sum_{i=1}^{r'}(a_i-\lfloor a_i/p\rfloor)$$
Pour conclure, consid\'erons $W\subset V$ une composante
irr\'eductible passant par $x$. Consid\'erons la stratification
naturelle
$$\Div_N(\P^1)=\cup_{\sum_i(a_i+1)=N+1}\sigma(a_1,\ldots,a_{r'})$$
avec $\sigma(a_1,\ldots,a_{r'})$ l'ensembles des diviseurs form\'es de
$r'$ points distincts avec les multiplicit\'es respectives
$a_1+1,\ldots,a_{r'}+1$ ($a_i\geq 0$). La cellule
$\overline{\sigma(a_1+1,\ldots,a_{r'}+1)}$ est irr\'eductible de dimension $r'$
; Supposons que
$$\Delta(W\times\PGL(1))\subset \overline{\sigma(a_1+1,\ldots,a_{r'}+1)}$$
et $\Delta(W)\cap\sigma(a_1+1,\ldots,a_{r'}+1)\not=\emptyset$. Comme la
fibre de $\Delta$ en un point g\'eom\'etrique de $\Delta(W)$ est de
dimension major\'ee par $\sum_{i=1}^{r'}(a_i-\lfloor a_i/p\rfloor)$, nous
obtenons
$$\dim W\times\PGL(1)\leq r'+ \sum_{i=1}^{r'}(a_i-\lfloor
a_i/p\rfloor)\leq\sum_{i=1}^{r'}(a_i+1)=N$$
Donc $\dim W\leq N-3$, ce qui termine la preuve de la proposition.
\subsection{Calcul de $\dim_k H^1(G,{\cal T}_C)$ et conclusion}
Revenons \`a la situation g\'en\'erale d'un rev\^etement cyclique de
degr\'e $p$ sur le corps $k$ : $C\rightarrow \Sigma=C/G$, avec
les points de branchement $b_1,\ldots,b_r$ et $m_1,\ldots,m_r$ leurs
conducteurs respectifs. Soit  $g_{\Sigma}$ le genre de la
courbe $\Sigma$.
\begin{lem}{\it
\index{formule g\'en\'erale de ramification} Soit $\gr\subset C$ le diviseur de
ramification, dont l'id\'eal en un point est la diff\'erente ; nous
avons la formule g\'en\'erale de ramification}
$$\Omega_C^1\cong {\cal O}_C(\gr)\otimes \pi^*(\Omega^1_{\Sigma})$$
\end{lem}
\begin{pr}
Le rev\^etement \'etant s\'eparable, nous avons la suite exacte
$$0\rightarrow\pi^*\Omega_{\Sigma}^1\rightarrow\Omega_C^1\rightarrow\Omega_{C/\Sigma}^1\rightarrow0$$
D'apr\`es \cite{Se1} III Proposition 14, ${\cal
  O}_C(\gr)\cong\Omega_C^1\otimes\pi^*(\Omega_{\Sigma}^1)^{-1}$. La
conclusion en d\'ecoule.
\end{pr}
\begin{prop}{\it Nous avons
$$\dim_kH^1(G,{\cal
  T}_C)=3g_{\Sigma}-3+\sum_{i=1}^r\Big\lceil{2\beta_i\over p}\Big\rceil$$
pour $\beta_i=(m_i+1)(p-1)$.}
\end{prop}
\begin{pr}
Soit $\gr\subset C$ le diviseur de ramification, ainsi d'apr\`es la
formule g\'en\'erale de ramification
$$\Omega_C^1\cong {\cal O}_C(\gr)\otimes \pi^*(\Omega^1_{\Sigma})$$
soit encore 
$${\cal T}_C\cong {\cal O}_C(-\gr)\otimes \pi^*({\cal T}_{\Sigma})$$
D'o\`u
$$\pi_*({\cal T}_C)\cong {\cal T}_{\Sigma}\otimes\pi_*({\cal
  O}_C(-\gr))\cong{\cal T}_{\Sigma}\otimes ({\cal O}_{\Sigma}\cap
  \pi_*({\cal O}_C(-\gr)))$$
Calculons d'abord le terme entre parenth\`ese ; soit $y\in \Sigma$. La
  fibre au point $y$ est\\
$${\cal O}_{\Sigma,y} \mbox{ si }
y\not\in\{b_1,\ldots,b_r\},\;\;\mbox{ et }\;\, {\cal
  O}_{\Sigma,b_i}\Big(\Big\lfloor {\beta_i\over p}\Big\rfloor b_i\Big)
\mbox{ si } y=b_i$$
car si $\gp(x_i)=b_i$, l'id\'eal de $\gr$ dans ${\cal O}_{C,x_i}$ est
la diff\'erente, donc est de valuation $\beta_i$. Donc
$$\pi_*^G({\cal T}_C)={\cal T}_{\Sigma}\otimes {\cal
  O}_{\Sigma}\Big(-\sum_{i=1}^r\Big\lfloor {\beta_i\over p}\Big\rfloor
  b_i\Big)$$
Ainsi la contribution ``globale'' \`a $H^1(G,{\cal T}_C)$ v\'erifie
$$H^1(\Sigma,\pi_*^G({\cal T}_C))\cong H^1\Big(\Sigma,{\cal
  O}_{\Sigma}\Big(\Big\lfloor {\beta_i\over p}\Big\rfloor
b_i\Big)\Big)\cong H^0\Big(\Sigma,\Omega_{\Sigma}^{\otimes2}\Big(-\sum_{i=1}^r\Big\lfloor {\beta_i\over p}\Big\rfloor
  b_i\Big)\Big)$$
le dernier isomorphisme provenant de la dualit\'e de Serre. Alors la
  formule de Riemann-Roch donne
$$\dim_kH^1(\Sigma,\pi_*^G({\cal
  T}_C))=3g_{\Sigma}-3+\sum_{i=1}^r\Big\lfloor {\beta_i\over
  p}\Big\rfloor$$
La contribution ``locale'' \`a $H^1(G,{\cal T}_C)$ est de dimension
$$\sum_{i=1}^r\dim_kH^1(G,\widehat{{\cal T}}_{C,x_i})=\sum_{i=1}^r
\Big\lceil{2\beta_i\over p}\Big\rceil-\Big\lfloor {\beta_i\over
  p}\Big\rfloor$$
D'o\`u
$$\dim_kH^1(G,{\cal
  T}_C)=3g_{\Sigma}-3+\sum_{i=1}^r\Big\lceil{2\beta_i\over
  p}\Big\rceil$$
\end{pr}
Nous pouvons maintenant obtenir la dimension de l'anneau de d\'eformation
versel local, associ\'e \`a un automorphisme $\sigma$ de conducteur $m$. Nous
supposons $m\geq 2$ et $p\geq 3$. Par le principe local-global
(Th\'eor\`eme \ref{thmloglo}), en supposant qu'il n'y a qu'un point de branchement, l'anneau de
d\'eformations universel global v\'erifie
$$\dim_{\Krull} R_{\glo}=N-2=\dim_{\Krull} R_{\sigma}+\dim_k
H^1(\Sigma,\pi_*^G({\cal T}_C))$$
D'apr\`es le calcul pr\'ec\'edent, nous obtenons ($g_{\Sigma}=0$) 
$$\dim_k H^1(\Sigma,\pi_*^G({\cal
  T}_C))=-3+\Big\lfloor{\beta\over p}\Big\rfloor$$
Comme $N+1=m+2$,
$$\dim_{\Krull}
R_{\sigma}=m+2-\Big\lfloor{\beta\over p}\Big\rfloor$$
D'o\`u le r\'esultat :
\begin{thm}
\label{thmGM}{\it
Soit $R_{\sigma}$ l'anneau de d\'eformations versel d\'efini par un
automorphisme $\sigma$ d'ordre $p$ et de conducteur $m=pq-l$ ($q\geq
1$ et $l\in [1,p-1]$). Nous avons
$$\dim_{\Krull} R_{\sigma}=\left\{\begin{array}{ll}q\mbox{ si }l\not=1\cr
q+1\mbox{ si } l=1\end{array}\right.$$
Et si $p>2$, $m<p-1$ et $(m,p)\not=(1,3)$, $R_{\sigma}$ est une intersection compl\`ete de dimension 1. }
\end{thm}
Il serait int\'eressant de savoir si la propri\'et\'e d'intersection
compl\`ete est toujours vraie pour d'autres valeurs de $m$ et de $p$
par exemple pour $m=p-1$.\\
Nous obtenons aussi le
corollaire suivant qui a pour cons\'equence un r\'esultat
de Green et Matignon (\cite{GrMa2}) lorsque $m<p-1$ (le r\'esultat de
Green et Matignon est encore valable pour $m=p-1$).

\begin{cor}{\it
Soit $S$ un anneau de valuation discr\`ete complet contenant $W(k)$. 
Supposons $m<p-1$. Alors il y a seulement un nombre fini de classes de
conjugaison d'automorphismes d'ordre $p$ de $S[[T]]$ qui induisent sur
la fibre sp\'eciale $k[[T]]$ un automorphisme de conducteur $m$.}
\end{cor}
\begin{pr}
Soit $\sigma$ un automorphisme de $k[[T]]$ de conducteur $m$.
La donn\'ee d'un automorphisme de
$S[[T]]$ qui induit $\sigma$, d\'efinit une d\'eformation de $\sigma$
et par cons\'equent un $W(k)$-automorphisme $$u:R_{\sigma}\rightarrow
S$$
La dimension de Krull de $R_{\sigma}$ est 1, donc $\dim_{\Krull}
R_{\sigma}/pR_{\sigma}=0$. Par cons\'equent $R_{\sigma}$ est un
  $W(k)$-morphisme de type fini. Supposons
  $R_{\sigma}=\sum_{i=1}^rW(k)e_i$. Notons $P_i(Y)$ un polyn\^ome
  unitaire \`a coefficients dans $W(k)$ tel que $P_i(e_i)=0$. Alors $u$
  est d\'etermin\'e par les images $y_i=u(e_i)\in S$ qui v\'erifient
  $P_i(y_i)=0$. Il n'y a qu'un nombre fini de solutions possibles pour
  un tel choix et donc pour $u.$
\end{pr} 
Il est possible de d\'eterminer dans certains cas les conducteurs
d'une d\'eformation g\'en\'erique :
\begin{cor}{\it
Soit $u$ et $v\in[1,p-1]$ d\'efinis par l'\'equation de B\'ezout
$up-vl=1$. Supposons que $v$ divise $p-1$ (par exemple si $l=1$).
Le rel\`evement de $(C,G)$ d\'efini par l'\'equation d'Artin-Schreier
d\'eform\'ee
$$
\xi^p-a(t)^{p-1}\xi=t^l,\;\;\;\; \mbox{ avec } a(t)=\left\{ \begin{array}{ll}
t^q+x_1t^{q-1}+\cdots+x_{q}&\mbox{ si } l=1\cr
t^q+x_1t^{q-1}+\cdots+x_{q-1}t&\mbox{ sinon }\cr
\end{array}
\right.$$
 a pour conducteur $q-1$ fois $p-1$ et une fois
$p-l$.} 
\end{cor}

\iffin
\end{document}